\documentclass[10pt]{article}
\usepackage{amssymb,latexsym}
\usepackage{epsfig}
\usepackage{eufrak}
\usepackage{amsmath}
\usepackage{mathrsfs}
\usepackage{color}

\newtheorem{theorem}{Theorem}
\newtheorem{lemma}{Lemma}

\newcommand{\be}{\begin{equation}}
\newcommand{\ee}{\end{equation}}
\newcommand{\bee}{\begin{eqnarray*}}
\newcommand{\eee}{\end{eqnarray*}}
\newcommand{\bel}{\begin{eqnarray}}
\newcommand{\eel}{\end{eqnarray}}
\newcommand{\bec}{\begin{cases}}
\newcommand{\eec}{\end{cases}}
\newcommand{\bem}{\begin{bmatrix}}
\newcommand{\eem}{\end{bmatrix}}

\newcommand{\la}{\label}
\newcommand{\li}{\left}
\newcommand{\ri}{\right}

\newcommand{\DEF}{\stackrel{\mathrm{def}}{=}}

\newcommand{\vep}{\varepsilon}
\newcommand{\lm}{\lambda}

\newcommand{\Up}{\Upsilon}

\newcommand{\si}{\sigma}

\newcommand{\de}{\delta}

\newcommand{\ga}{\gamma}

\newcommand{\vse}{\vartheta}
\newcommand{\se}{\theta}
\newcommand{\Se}{\Theta}

\newcommand{\ze}{\zeta}
\newcommand{\al}{\alpha}
\newcommand{\ba}{\beta}

\newcommand{\ro}{\rho}
\newcommand{\ka}{\kappa}
\newcommand{\om}{\omega}
\newcommand{\Om}{\Omega}

\newcommand{\f}{\frac}
\newcommand{\sq}{\sqrt}
\newcommand{\cd}{\cdots}

\newcommand{\qu}{\quad}
\newcommand{\qqu}{\qquad}
\newcommand{\fa}{\forall}

\newcommand{\mscr}{\mathscr}
\newcommand{\mcal}{\mathcal}
\newcommand{\mbf}{\mathbf}
\newcommand{\bb}{\mathbb}

\newcommand{\wh}{\widehat}
\newcommand{\wt}{\widetilde}
\newcommand{\mrm}{\mathrm}
\newcommand{\bs}{\boldsymbol}

\newcommand{\tx}{\text}

\newcommand{\iy}{\infty}

\newcommand{\pa}{\partial}

\newcommand{\bed}{\begin{description}}
\newcommand{\eed}{\end{description}}
\newcommand{\bei}{\begin{itemize}}
\newcommand{\eei}{\end{itemize}}
\newcommand{\ben}{\begin{enumerate}}
\newcommand{\een}{\end{enumerate}}
\newcommand{\bib}{\bibitem}
\newcommand{\beL}{\begin{lemma}}
\newcommand{\eeL}{\end{lemma}}
\newcommand{\beT}{\begin{theorem}}
\newcommand{\eeT}{\end{theorem}}
\newcommand{\sect}{\section}

\newcommand{\bpf}{\begin{pf}}
\newcommand{\epf}{\end{pf}}

\newcommand{\bi}{\binom}

\newcommand{\pfbox}{\hfill\mbox{$\Box$}}

\newenvironment{pf}{\paragraph*{Proof{\rm.}}}{\pfbox\bigskip}

\begin{document}

\title{{\bf Consecutive Sequential Probability Ratio Tests of Multiple Statistical Hypotheses} \thanks{The author had been previously working with Louisiana
State University at Baton Rouge, LA 70803, USA, and is now with Department of Electrical Engineering, Southern University and A\&M College,
Baton Rouge, LA 70813, USA; Email: chenxinjia@gmail.com. The main results of this paper have appeared in Proceedings of SPIE Conferences,
Baltimore, Maryland,  April 24-27, 2012.} }

\author{Xinjia Chen}

\date{First submitted in June 2012}

\maketitle

\begin{abstract}

In this paper, we develop a simple approach for testing multiple statistical hypotheses based on the observations of a number of probability
ratios enumerated consecutively with respect to the index of hypotheses. Explicit and tight bounds for the probability of making wrong decisions
are obtained for choosing appropriate parameters for the proposed tests.  In the special case of testing two hypotheses, our tests reduce to
Wald's sequential probability ratio tests.

\end{abstract}


\section{Introduction}

Consider a continuous-time stochastic process $(X_t)_{t \in [0, \iy) }$ defined in a probability space $(\Om, \mscr{F}, \Pr )$. Suppose that the
stochastic process $(X_t)_{t \in [0, \iy) }$ is parameterized by $\se \in \Se$.  In many applications of engineering and sciences, it is
desirable to infer the true value of $\se$ based on the observation of such stochastic processes.  This topic can be formulated as a general
problem of testing $m$ mutually exclusive and exhaustive composite hypotheses: \be \la{mainpr} \mscr{H}_0: \se \in {\Se}_0, \qu \mscr{H}_1: \se
\in {\Se}_1, \qu \ldots, \qu \mscr{H}_{m - 1}: \se \in {\Se}_{m - 1}, \ee where $\Se_i = \{ \se \in \Se: \se_i < \se \leq \se_{i+1} \}, \; i =
0, 1, \cd, m - 1$ with $- \iy = \se_0 < \se_1 < \cd < \se_{m - 1} < \se_m = \iy$. To control the probabilities of making wrong decisions, for
pre-specified numbers $\de_i \in (0, 1), \; i = 0, 1, \cd, m - 1$, it is typically required that \be \la{mainreq}
 \Pr \{ \tx{Reject} \; \mscr{H}_i \mid \se \} \leq \de_i, \qqu \fa \se
\in \varTheta_i, \qu i = 0, 1, \cd, m - 1 \ee  where $\varTheta_i = \{ \se \in \Se_i: \se_i^{\prime \prime} \leq \se \leq \se_{i+1}^{\prime} \},
\; i = 0, 1, \cd, m - 1$ with $\se_i^\prime, \se_i^{\prime \prime} \in \Se, \; i = 1, \cd, m - 1$ satisfying $- \iy = \se_0^{\prime \prime} <
\se_i^\prime < \se_i < \se_i^{\prime \prime} \leq \se_{i + 1}^\prime < \se_{i + 1} < \se_{i+1}^{\prime \prime} < \se_m^\prime = \iy$ for $i = 1,
\cd, m - 2$. The set $\cup_{i = 1}^{m-1} (\se_i^\prime, \se_i^{\prime \prime})$ is referred to as the indifference zone, since no specification
on risk is imposed for the set. Here we consider continuous-time processes for the sake of generality, since discrete-time stochastic processes
can be treated as right-continuous processes in continuous time.

The hypothesis testing problem defined by (\ref{mainpr}) and (\ref{mainreq}) has been studied extensively for more than a half century (see,
\cite{Eisenberg, Ghosh} and the references therein). In particular, for the special problem of testing two hypotheses, Wald \cite{Wald} invented
the famous Sequential Probability Ratio Tests (SPRTs). Armitage \cite{Armitage} extended Wald's SPRTs to the general problem of testing multiple
hypotheses. Lorden \cite{Lorden} proposed sequential likelihood ratio tests for the same problem.  Baum \cite{Baum} established  multiple
sequential probability ratio tests in a Bayesian framework.  At present the general theory of tests on multiple statistical hypotheses is much
less developed than for the two-decision situation.  Existing methods suffer from one or more of the following drawbacks: (i) There is no
rigorous method for controlling the risk of making wrong decisions; (ii) The method of bounding the risk of making wrong decisions is too
conservative; (iii) The application is limited to simple hypotheses; (iv) The application is limited by the number of hypotheses.   Motivated by
this situation, we develop a new class of tests, referred to as {\it Consecutive Sequential Probability Ratio Tests} (CSPRTs) based on the
principle of probabilistic comparison proposed in \cite{Chenlimits, Chenhypothesis, ChenSPIE2012}.

The remainder of this paper is organized as follows.  In Section \ref{secmul}, we introduce the connection between multi-hypotheses testing and
sequential random intervals.  In Section \ref{secPri}, we describe the principle of probabilistic comparison. In Section \ref{secCSPRT}, we
apply the principle of probabilistic comparison to develop consecutive sequential probability ratio tests.   In Section \ref{secCon},  we
establish consecutive sequential probability ratio tests on parameters of continuous-time processes.  Section \ref{secCons} is the conclusion.
All proofs are given in Appendices. The main results of this paper have been appeared in our conference paper \cite{ChenSPIE2012}.

Throughout this paper, we shall use the following notations. The empty set is denoted by $\emptyset$. The set of positive integers is denoted by
$\bb{N}$.  The notation $\Pr \{ E \mid \se \}$ denotes the probability of the event $E$ associated with parameter $\se$.  The expectation of a
random variable is denoted by $\bb{E}[.]$. The support of a random variable $Z$ is denoted by $I_Z$.  In the discrete-time case, the stochastic
process $(X_t)_{t \in [0, \iy) }$ is actually a sequence of random variables $X_1, X_2, \cd$.  For simplicity of notations, let $\bs{\mcal{X}}_n
= (X_1, \cd, X_n)$ for $n \in \bb{N}$. Let $\mbf{x}_n = (x_1, \cd, x_n)$ denote the realization of $\bs{\mcal{X}}_n$. Let $f_n (\mbf{x}_n; \se)$
denote the probability density function (PDF) or probability mass function (PMF) of $(X_1, \cd, X_n)$ parameterized by $\se \in \Se$.
Accordingly, replacing $\mbf{x}_n$ in $f_n (\mbf{x}_n; \se)$ by $\bs{\mcal{X}}_n$ gives the likelihood function $f_n (\bs{\mcal{X}}_n; \se)$.
For $\se^{\prime}, \se^{\prime \prime} \in \Se$ and $\ka > 0$, we use $\Up_n (\bs{\mcal{X}}_n; \; \se^{\prime}, \se^{\prime \prime} ) \bs{\sim}
\ka$  to represent $f_n (\bs{\mcal{X}}_n; \; \se^{\prime \prime}) \bs{\sim} \ka f_n (\bs{\mcal{X}}_n; \; \se^{\prime} )$, where  ``$\bs{\sim}$''
is a relation such as ``$<, \; = , \;
>, \; \leq, \; \geq$'', corresponding to ``less than, equal, greater than, less or equal, greater or equal'', respectively.
The notation $\Up_n (\bs{\mcal{X}}_n; \; \se^{\prime}, \se^{\prime \prime} )$ can be interpreted as the likelihood ratio $\f{ f_n
(\bs{\mcal{X}}_n; \; \se^{\prime \prime}) }  { f_n (\bs{\mcal{X}}_n; \; \se^{\prime} )  }$ whenever $f_n (\bs{\mcal{X}}_n; \; \se^{\prime} )$ is
not equal to $0$.  We shall frequently use the concept of unimodal function.  A function is said to be unimodal with respect to $\se \in \Se$ if
there exists a number $\se^*$ such that the function is non-decreasing with respect to $\se \in \Se$ no greater than $\se^*$ and is
non-increasing with respect to $\se \in \Se$ no less than $\se^*$.  The other notations and concepts will be made clear as we proceed.

\section{Multi-hypotheses Testing and Sequential Random Intervals}  \la{secmul}

As demonstrated in \cite{Chenlimits}, the general hypothesis testing problem defined by (\ref{mainpr}) and (\ref{mainreq}) can be cast into the
framework of constructing a sequential random interval with pre-specified coverage probabilities. This can be illustrated in the sequel.

To reach a fast decision, it is desirable to solve the hypothesis testing problem by a multistage approach such that the sampling procedure is
divided into $s$ stages with observational times $t_{\ell}, \; \ell = 1, \cd, s$, where $t_{\ell}$ is the observational time at the $\ell$-th
stage. Starting from $\ell = 1$, at the $\ell$-th stage, based on the observation of $(X_t)_{0 \leq t \leq t_{\ell} }$, pre-determined stopping
and decision rules are applied to check whether the accumulated observational data is sufficient to accept a hypothesis and terminate the
sampling procedure. If the observational data is considered to be insufficient for making a decision, then proceed to the next stage of
observation.  The observation is continued stage by stage until a hypothesis is accepted at some stage. Although the number of stages $s$ may be
infinity, for practical considerations, the stopping and decision rules are required to guarantee that the sampling procedure will surely
eventually terminate with a finite number of stages. Central to a multistage procedure are the stopping and decision rules, which can be related
to a sequential random interval described as follows.  Let $\se_0^\prime = -\iy$ and $\se_m^{\prime \prime} = \iy$. For $i = 0, 1, \cd, m - 1$,
let $\mscr{I}_i$ denote the open interval {\small $( \se_i^\prime, \se_{i+1}^{\prime \prime} )$}. Let $\bs{l}$ be the index of stage at the
termination of the sampling procedure. Let $\bs{\mcal{L}}$ and $\bs{\mcal{U}}$ be random variables defined in terms of samples of the stochastic
process up to the $\bs{l}$-th stage such that the sequential random interval {\small $( \bs{\mcal{L}}, \bs{\mcal{U}})$} has $m$ possible
outcomes $\mscr{I}_i, \; i = 0, 1, \cd, m - 1$ and that $\Pr \{ \bs{\mcal{L}}  < \se < \bs{\mcal{U}}  \mid \se \} > 1 - \de_i$ for any $\se \in
\varTheta_i$ and $i = 0, 1, \cd, m - 1$.  Given that the sequential random interval $( \bs{\mcal{L}}, \bs{\mcal{U}} )$ satisfying such
requirements is constructed, the risk requirement (\ref{mainreq}) can be satisfied by using $( \bs{\mcal{L}}, \bs{\mcal{U}} )$ to define a
decision rule such that, for $i = 0, 1, \cd, m - 1$, hypothesis $\mscr{H}_i$ is accepted when the sequential random interval $( \bs{\mcal{L}},
\bs{\mcal{U}} )$ takes $\mscr{I}_i$ as its outcome at the termination of the sampling process. It follows that $\{ \tx{Accept} \; \mscr{H}_i \}
=\{ \bs{\mcal{L}} < \se < \bs{\mcal{U}} \}$ for any $\se \in \varTheta_i$ and $i = 0, 1, \cd, m - 1$. Therefore, to solve the multi-valued
decision problem defined by (\ref{mainpr}) and (\ref{mainreq}), the objective is to ensure that $\se$ is included in the sequential random
interval with pre-specified probabilities.  In the sequel, we shall propose a general approach for defining stopping and decision rules for the
construction of such sequential random interval.

\section{Principle of Probabilistic Comparison} \la{secPri}

In \cite{Chenlimits, Chenhypothesis, ChenSPIE2012}, a general methodology has been proposed for constructing sequential random intervals with
prescribed specifications of coverage probabilities.  The main idea is to use one-sided confidence sequences to control the coverage probability
of the sequential random interval.  Assume that the number of stages, $s$, and the observational times, $t_{\ell}, \; \ell = 1, \cd, s$, are
given. Assume that for $\ell = 1, \cd, s$ and $i = 1, \cd, m-1$, random variables $L_{\ell, i}$ and $U_{\ell, i}$ can be defined in terms of
positive numbers $\ze, \; \al_i, \; \ba_i$ and the set of random variables $(X_t)_{0 \leq t \leq t_{\ell} }$ such that $\Pr \{ L_{\ell, i} \geq
\se \mid \se \}$ and $\Pr \{ U_{\ell, i} \leq \se \mid \se \}$  can be made arbitrarily small by decreasing $\ze \al_i$ and $\ze \ba_i$
respectively. Due to such assumption, we call $(-\iy, L_{\ell, i}]$ and $[U_{\ell, i}, \iy)$ one-sided confidence intervals for $\se$.
Accordingly, $(-\iy, L_{\ell, i}], \; \ell = 1, \cd, s$ and $[U_{\ell, i}, \iy), \; \ell = 1, \cd, s$ are said to be one-sided confidence
sequences for $\se$. In view of the controllability of the coverage probabilities of the one-sided confidence intervals, the number $\ze$ is
referred to as the {\it coverage tuning parameter}, and $\al_i, \; \ba_i, \; i = 1, \cd, m-1$ are called {\it weighting coefficients}.  Given
that $\ze$ is sufficiently small, $\se
> \se_i^\prime$ will be credible if $L_{\ell, i} > \se_i^\prime$ is observed. Similarly, $\se < \se_i^{\prime \prime}$ will be credible if
$U_{\ell, i} < \se_i^{\prime \prime}$ is observed.  To figure out the general structure of stopping and decision rules, imagine that the
sampling procedure is stopped at the $\ell$-th stage and $\mscr{I}_i$ is to be designated as the outcome of the sequential random interval.
Since $\mscr{I}_i$ contains $[\se_i^{\prime \prime}, \se_{i + 1}^{\prime}]$, it follows that for $\se \in [\se_i^{\prime \prime}, \se_{i +
1}^{\prime}]$, it is true that $\se < \se_j^{\prime \prime}$ for $j > i$ and $\se
> \se_j^\prime$ for $j \leq i$.  This implies that, if the coverage tuning parameter $\ze$ is sufficiently small, then it is very likely to observe
that $U_{\ell, j} < \se_j^{\prime \prime}$ for $j > i$ and $L_{\ell, j} > \se_j^\prime$ for $j \leq i$. Therefore, turning this thinking around
leads to the following stopping and decision rules:

\vspace{0.05in}

\begin{tabular} {|l |}
\hline $ \tx{ {\it Continue observing the stochastic processes until for some $i \in \{0, 1, \cd, m - 1\}$, the event}}$\\
$\tx{ {\it $\{ U_{\ell, j} < \se_j^{\prime \prime} \; \tx{for} \;  j > i \;  \tx{and} \; L_{\ell, j} > \se_j^\prime \;  \tx{for} \;  j \leq i
\}$ occurs at some stage with index $\ell \in \{1, \cd, s\}$}}$.\\
$\tx{ {\it At the termination of the sampling process, make the following decision: If such index $i$ is unique,}}$\\
$ \tx{ {\it  then designate $\mscr{I}_i$ as the outcome of the sequential random interval. If there are multiple indexes  }}$\\
$ \tx{ {\it  satisfying the condition, then pick one of them and assign the corresponding interval $\mscr{I}_i$ as the}}$\\
$ \tx{ {\it   outcome of the sequential random interval based on a predetermined policy.}}$
\\ \hline
\end{tabular}

\vspace{0.05in}

 The idea in the derivation of the above stopping and decision rules is to infer
 the location of $\se$ relative to the sequential random interval by comparing
 the confidence limits with the endpoints of the sequential random interval. Due to
the probabilistic nature of the comparison, such method of constructing stopping and decision rules is referred to as the {\it Principle of
Probabilistic Comparison}.  It should be noted that similar principles have been proposed in \cite{Chenestimation, Chen_rule} for multistage
estimation of parameters.   The properties of the above  stopping and decision rules are indicated by the following probabilistic result.

\beT \la{Multi-Valued Inclusion Principle}

Let $a_0 = b_0 = - \iy, \; a_m = b_m = \iy$ and  $a_i < b_i \leq a_{i + 1} < b_{i + 1}$ for $i = 1, \cd, m - 2$.    Let $\varTheta_0 = (- \iy,
a_1], \; \varTheta_{m-1} = [b_{m - 1}, \iy )$ and $\varTheta_i = [b_i, a_{i + 1}]$ for $i = 1, \cd, m - 2$.  Let $(\Om, \mscr{F}, \{
\mscr{F}_\ell \}, \Pr )$ be a filtered space. Let $\bs{\tau}$ be a proper stopping time with a support $I_{\bs{\tau}}$.  For $\ell \in
I_{\bs{\tau}}$, let $L_{\ell, m } = - \iy, \; U_{\ell, 0}  = \iy$ and let $L_{\ell, i}, \; U_{\ell, i}, \; i = 1,\cd, m - 1$ be random variables
measurable in $\mscr{F}_\ell$.  Let $\bs{\mcal{L}}$ and $\bs{\mcal{U}}$ be random variables such that $\cup_{i = 0}^{m-1} \{ \bs{\mcal{L}} =
a_i, \; \bs{\mcal{U}} = b_{i + 1} \} = \Om$ and that {\small $\{ \bs{\tau} = \ell, \; \bs{\mcal{L}} = a_j, \; \bs{\mcal{U}} = b_{j + 1} \}
\subseteq \{ L_{\ell, i} \geq a_i, \; 0 < i \leq j \; \tx{and} \; U_{\ell, i} \leq b_i, \; j < i < m \}$ } for $\ell \in I_{\bs{\tau}}$ and $j =
0, 1, \cd, m - 1$. Then,  $\Pr \{ \bs{\mcal{L}} \geq \se \} = \Pr \{ \bs{\mcal{L}} \geq a_{i+1} \} \leq \Pr \{ L_{\ell, i + 1} \geq a_{i+1} \;
\tx{for some} \; \ell \in I_{\bs{\tau}} \}$ and $\Pr \{ \bs{\mcal{U}} \leq \se \} = \Pr \{ \bs{\mcal{U}} \leq b_i \} \leq \Pr \{ U_{\ell, i}
\leq b_i \; \tx{for some} \; \ell \in I_{\bs{\tau}} \}$ for $i = 0, 1, \cd, m - 1$ and $\se \in \varTheta_i$.

\eeT

See Appendix \ref{Multi-Valued Inclusion Principle_app} for a proof.

\section{Consecutive Sequential Probability Ratio Tests}  \la{secCSPRT}

In this section, we shall apply the principle of probabilistic comparison and Theorem \ref{Multi-Valued Inclusion Principle} to develop a new
class of tests for solving the multi-valued decision problem defined by (\ref{mainpr}) and (\ref{mainreq}) regarding the parameter $\se \in \Se$
associated with a discrete process $( X_n )_{n \in \mcal{N}}$, where $\mcal{N}$ is a subset of positive integers. For generality,  we do not
restrict $\mcal{N}$ as an unbounded set such as $\bb{N}$.  Our purpose is to accommodate the situation that the sequence of $X_n$ can be of
finite length.  A familiar example can be found in the context of sampling without replacement from a finite population of $N$ units, among
which $N p$ units having a certain attribute.  If we define a Bernoulli random variable $X_n$ such that $X_n$ assumes values $1$  or $0$ in
accordance with whether the $n$-th drawn unit has the attribute, then we have a sequence of dependent Bernoulli random variables $(X_n)_{n \in
\mcal{N}}$ with $\mcal{N} = \{ 1, 2, \cd, N \}$. Throughout the remainder of this paper, we use symbol $N^*$ to denote $\iy$ if $\mcal{N}$ is
unbounded and otherwise the maximum of $\mcal{N}$.

\subsection{Confidence Sequences}

For the purpose of deriving sequential tests based on the principle of probabilistic comparison, we need a method for constructing confidence
sequences as described by the following theorem.

\beT  \la{CI sequence} For $n \in \mcal{N}$, let $\bs{\mcal{X}}_n$ be random variables parameterized by $\se \in \Se$ and let the likelihood
function be denoted by $f_n(\bs{\mcal{X}}_n; \se)$.  Let $\de \in (0, 1)$ and let $\se_0, \; \se_1 \in \Se$ with $\se_0 < \se_1$. Define random
variables $L_n (\bs{\mcal{X}}_n) = \inf  \{ \vse \in \Se: \Up_n (\bs{\mcal{X}}_n; \se_1,  \vse) > \f{\de}{2} \}$ and $U_n (\bs{\mcal{X}}_n) =
\sup \{ \vse \in \Se:
 \Up_n (\bs{\mcal{X}}_n;  \se_0,  \vse) > \f{\de}{2} \}$.
The following statements hold true.

(I) For all $\se \in \Se$,  \bee &  & \Pr \{ L_n (\bs{\mcal{X}}_n) \leq \se  \; \tx{for all} \; n \in \mcal{N} \mid \se \} \geq 1 - \f{\de}{2}, \\
&  & \Pr \{  U_n (\bs{\mcal{X}}_n) \geq \se \; \tx{for all} \; n  \in \mcal{N} \mid \se \} \geq 1 - \f{\de}{2}, \\
&  & \Pr \{ L_n (\bs{\mcal{X}}_n) \leq \se \leq U_n (\bs{\mcal{X}}_n) \; \tx{for all} \; n  \in \mcal{N} \mid \se \} \geq 1 - \de. \eee

(II) For all $n \in \mcal{N}$,
\[
 \{ L_n (\bs{\mcal{X}}_n) > \se_0 \} \subseteq \li \{ \Up_n (\bs{\mcal{X}}_n; \se_1,  \se_0) \leq \f{\de}{2} \ri \}, \qqu
 \{ U_n (\bs{\mcal{X}}_n) < \se_1 \} \subseteq \li \{ \Up_n (\bs{\mcal{X}}_n; \se_0,  \se_1) \leq \f{\de}{2} \ri \}.
\]

(III) If $f_n (\bs{\mcal{X}}_n; \se)$ is unimodal with respect to $\se \in \Se$, then
\[
 \{ L_n (\bs{\mcal{X}}_n) \geq \se_0 \} \supseteq \li \{ \Up_n (\bs{\mcal{X}}_n; \se_1,  \se_0) \leq \f{\de}{2} \ri \}, \qqu
 \{ U_n (\bs{\mcal{X}}_n) \leq \se_1 \} \supseteq \li \{ \Up_n (\bs{\mcal{X}}_n; \se_0,  \se_1) \leq \f{\de}{2} \ri \}
\]
for all $n \in \mcal{N}$.  \eeT

See Appendix \ref{CI sequence_app} for a proof.

Assuming  that $f_n (\bs{\mcal{X}}_n; \se)$ is unimodal with respect to $\se \in \Se$ and that $\Se$ is a discrete set or $f_n (\bs{\mcal{X}}_n;
\se)$ is continuous with respect to $\se \in \Se$, we have
\[
 \{ L_n (\bs{\mcal{X}}_n) \geq \se_0 \} = \li \{ \Up_n (\bs{\mcal{X}}_n; \se_1,  \se_0) \leq \f{\de}{2} \ri \}, \qqu
 \{ U_n (\bs{\mcal{X}}_n) \leq \se_1 \} = \li \{ \Up_n (\bs{\mcal{X}}_n; \se_0,  \se_1) \leq \f{\de}{2} \ri \}
\]
for all $n \in \mcal{N}$.

\subsection{CSPRTs on Multiple Composite Hypotheses} \la{CSPRTSC}

In order to construct a sequential test,  choose  $\al_i, \; \ba_i \in (0, 1)$ for $i = 1, \cd, m - 1$ and $\al_m = \ba_0 = 0$. Define lower
confidence limit
\[ L_{n, i} = \inf \{ \vse \in \Se: \Up_n (\bs{\mcal{X}}_n; \se_i^{\prime \prime},  \vse) > \al_i \}
\]
and upper confidence limit
\[
U_{n, i} = \sup \{ \vse \in \Se: \Up_n (\bs{\mcal{X}}_n;  \se_i^{\prime},  \vse) > \ba_i \}
\]
for $i = 1, \cd, m - 1$.  Making use of the principle of probabilistic comparison, we propose stopping and decision rules as follows:

Continue the sampling process until there exists an index $j \in \{0, 1, \cd, m - 1 \}$ such that
\[
L_{n, i} \geq \se_i^{\prime}, \qu 0 <  i \leq j \qu \tx{and} \qu U_{n, i} \leq \se_i^{\prime \prime}, \qu j < i < m.
\]
At the termination of the sampling process, accept $\mscr{H}_j$ with the index $j$ satisfying the stopping condition.

As a consequence of Theorems \ref{Multi-Valued Inclusion Principle} and \ref{CI sequence}, we have that if the sampling process will eventually
terminate with probability $1$, then  $\Pr \{ \tx{Reject} \; \mscr{H}_i \mid \se \} \leq \al_{i + 1} + \ba_i$ for $0 \leq i < m$ and $\se \in
\varTheta_i$.

Under the assumption that $f_n (\bs{\mcal{X}}_n; \se)$ is unimodal with respect to $\se \in \Se$ and that $\Se$ is a discrete set or $f_n
(\bs{\mcal{X}}_n; \se)$ is continuous with respect to $\se \in \Se$, it follows from Theorem \ref{CI sequence} that $\{ L_{n, i} \geq
\se_i^{\prime} \} = \{ \Up_n (\bs{\mcal{X}}_n; \se_i^{\prime \prime}, \se_i^{\prime}) \leq \al_i \}$ and $\{ U_{n, i} \leq \se_i^{\prime \prime}
\} = \{  \Up_n (\bs{\mcal{X}}_n; \se_i^{\prime}, \se_i^{\prime \prime}) \leq \ba_i \}$ for $0 < i < m$. Hence, the stopping and decision rules
can be simplified as follows:

\vspace{0.05in}

\begin{tabular} {|l |}
\hline $ \tx{ {\it Continue the sampling process until there exists an index $j$ in the set $\{0, 1, \cd, m - 1 \}$ such that}}$\\
$\tx{ {\it $\Up_n (\bs{\mcal{X}}_n; \se_i^{\prime}, \se_i^{\prime \prime}) \geq \f{1}{\al_i}$ for $0 < i \leq j$ and $\Up_n (\bs{\mcal{X}}_n;
\se_i^{\prime}, \se_i^{\prime \prime}) \leq \ba_i$ for $j < i < m$.
At the termination of}}$\\
$\tx{ {\it the sampling process, accept $\mscr{H}_j$ with the index $j$ satisfying the stopping condition}}$.
\\ \hline
\end{tabular}

\vspace{0.05in}

A salient feature of our test is that $m - 1$ consecutive probability ratios are used for defining the stopping and decision rules.  The name
{\it Consecutive Sequential Probability Ratio Test} is derived from such nature of the test.  We have established that the consecutive
sequential probability ratio test has the following properties.

\beT \la{MSPRT_Composite}

Assume that the likelihood function $f_n (\bs{\mcal{X}}_n; \se)$ is unimodal with respect to $\se \in \Se$ for any $n \in \mcal{N}$. If the
sampling process will eventually terminate according to the stopping rule with probability $1$, then the following statements (I)--(III) hold
true:

(I) $\Pr \{ \tx{Reject} \; \mscr{H}_i \mid \se \} \leq \al_{i + 1} + \ba_i$ for $0 \leq i \leq m - 1$ and $\se \in \varTheta_i$.

(II) For $j = 1, \cd, m - 1$, $\Pr \{ \tx{Accept} \; \mscr{H}_i \; \tx{with some index $i$ no less than $j$} \mid \se \}$ is no greater than
$\al_j$ and is non-decreasing with respect to $\se \in \Se$ no greater than $\se_j^\prime$.

(III) For $j = 1, \cd, m - 1$, $\Pr \{ \tx{Accept} \; \mscr{H}_i \; \tx{with some index $i$ less than $j$} \mid \se \}$ is no greater than
$\ba_j$ and is non-increasing with respect to $\se \in \Se$ no less than $\se_j^{\prime \prime}$.

Moreover, the sampling process will eventually terminate according to the stopping rule with  probability $1$,  provided that the following
additional assumption is satisfied:
 For arbitrary $\al, \ba \in (0, 1)$ and $\se \in \Se$, \be \la{TwoStop}
 \Pr \li \{ \ba < \Up_n (\bs{\mcal{X}}_n; \se_i^{\prime}, \se_i^{\prime \prime})  < \f{1}{\al}  \mid \se \ri \}  \to 0, \qqu i = 1, \cd, m - 1
\ee as the sample number $n$ tends to $N^*$.

 \eeT

 See Appendix \ref{MSPRT_Composite_app} for a proof.  It should be emphasized that throughout this paper,  the notion of ``the
sampling process will eventually terminate according to the stopping rule'' is that the stopping rule is satisfied for some $n \in \mcal{N}$.

Statement (I) of Theorem \ref{MSPRT_Composite} provides a simple method for controlling the risk of making wrong decisions.  To satisfy the risk
requirement (\ref{mainreq}), it suffices to choose $\al_i$ and $\ba_i$ such that $\al_{i + 1} + \ba_i \leq \de_i$ for $0 \leq i \leq m - 1$.
Specially, one can simply use $\al_1 = \de_0, \; \ba_{m-1} = \de_{m -1}$ and $\al_{i+1} = \ba_{i} = \f{\de_i}{2}$ for $1 \leq i \leq m - 2$ in
the stopping and decision rules for purpose of ensuring (\ref{mainreq}).

\subsection{CSPRTs on Multiple Simple Hypotheses} \la{CSPRTSS}

 In some situations, it may be interesting to test multiple simple hypotheses
 \be
 \la{simdefA}
 \mscr{H}_0 : \se =
\se_0, \qqu \mscr{H}_1 : \se = \se_1, \qqu \cd, \qqu \mscr{H}_{m-1} : \se = \se_{m-1}. \ee
 For risk control purpose, it is typically required that, for prescribed numbers $\de_i \in (0, 1)$,
 \be
 \la{simdefB}
 \Pr
\li \{ \tx{Reject} \; \mscr{H}_i \mid \se_i \ri \} \leq \de_i, \qqu  i = 0, 1, \cd, m - 1. \ee
 As before, let $\al_i, \ba_i \in (0, 1)$ for $i = 1,
\cd, m - 1$ and $\al_{m} = \ba_0 = 0$. Define lower confidence limit \[ L_{n, i} = \inf  \{ \vse \in \Se: \Up_n (\bs{\mcal{X}}_n; \se_{i + 1},
\vse)
> \al_i \}
\]
and upper confidence limit
\[
U_{n, i} = \sup \{ \vse \in \Se: \Up_n (\bs{\mcal{X}}_n;  \se_i,  \vse) > \ba_i \}
\]
for $i = 1, \cd, m - 1$.  By the principle of probabilistic comparison, we propose the following stopping and decision rules:

Continue the sampling process until  there exists an index $j \in \{0, 1, \cd, m - 1 \}$ such that
\[
L_{n, i} \geq \se_i, \qu 0 \leq i < j \qu \tx{and} \qu U_{n, i} \leq \se_{i + 1}, \qu j \leq i \leq m - 2.
\]
At the termination of the sampling process, accept $\mscr{H}_j$ with the index $j$ satisfying the stopping condition.

Under the assumption that $f_n (\bs{\mcal{X}}_n; \se)$ is unimodal with respect to $\se \in \Se$, it follows from Theorem \ref{CI sequence} that
$\{ L_{n, i} \geq \se_i  \} = \{ \Up_n (\bs{\mcal{X}}_n; \se_{i + 1}, \se_i) \leq \al_i \}$ and $\{ U_{n, i} \leq \se_{i + 1} \} = \{ \Up_n
(\bs{\mcal{X}}_n; \se_i, \se_{i + 1}) \leq \ba_i \}$ for $0 \leq i \leq m - 2$. Hence, the stopping and decision rules can be simplified  as
follows:

\vspace{0.05in}

\begin{tabular} {|l |}
\hline $ \tx{ {\it Continue the sampling process until there exists an index $j$ in the set $\{0, 1, \cd, m - 1 \}$ such that}}$\\
$\tx{ {\it $\Up_n (\bs{\mcal{X}}_n; \se_{i-1}, \se_{i}) \geq \f{1}{\al_i}$ for $0 < i \leq j$ and $\Up_n (\bs{\mcal{X}}_n; \se_{i-1}, \se_{i})
\leq \ba_i$ for $j < i < m$.
At the termination}}$\\
$\tx{ {\it of the sampling process, accept $\mscr{H}_j$ with the index $j$ satisfying the stopping condition}}$.
\\ \hline
\end{tabular}

\vspace{0.05in}

We have shown that the above consecutive sequential probability ratio test has the following properties.

\beT  \la{MSPRT_Simple}  If the sampling process will eventually terminate according to the stopping rule with probability $1$, then $\Pr \{
\tx{Reject} \; \mscr{H}_i \mid \se_i \} \leq \al_{i+1} + \ba_i$ for $0 \leq i \leq m - 1$.  Moreover, the sampling process will eventually
terminate according to the stopping rule with probability $1$, provided that the likelihood function $f_n (\bs{\mcal{X}}_n; \se)$ is unimodal
with respect to $\se \in \Se$ for any positive integer $n$, and that for arbitrary $\al, \ba \in (0, 1)$ and $\se \in \Se$,
\[
\Pr \li \{ \ba < \Up_n (\bs{\mcal{X}}_n; \se_{i-1}, \se_{i}) < \f{1}{\al} \mid \se \ri \} \to 1, \qqu i = 1 \cd, m - 1 \] as the sample number
$n$ tends to $N^*$.  \eeT

See Appendix \ref{MSPRT_Simple_app} for a proof.

According to Theorem \ref{MSPRT_Simple},   to guarantee the risk requirement (\ref{mainreq}), it suffices to choose $\al_i$ and $\ba_i$ such
that $\al_{i + 1} + \ba_i \leq \de_i$ for $0 \leq i \leq m - 1$.   Particularly, one can use $\al_1 = \de_0, \; \ba_{m-1} = \de_{m -1}$ and
$\al_{i+1} = \ba_{i} = \f{\de_i}{2}$ for $1 \leq i \leq m - 2$ in the stopping and decision rules to ensure that $\Pr \li \{ \tx{Reject} \;
\mscr{H}_i \mid \se_i \ri \} \leq \de_i, \; i = 0, 1, \cd, m - 1$.

\subsection{General Termination Properties}

In Theorems \ref{MSPRT_Composite} and \ref{MSPRT_Simple}, one of the assumptions that we use to establish the termination properties is that the
likelihood functions are unimodal on $\Se$.  Actually, with regard to the CSPRTs on composite and simple hypotheses proposed in Sections
\ref{CSPRTSC} and \ref{CSPRTSS}, the termination properties are valid under fairy general assumptions, as asserted by the following results.

\beT
 \la{genSTgen}

The sampling process will eventually terminate according to the stopping rule with probability $1$, provided that the following assumptions are
satisfied:

(I) For arbitrary $\al, \ba \in (0, 1)$ and $\se, \se^{\prime}, \se^{\prime \prime} \in \Se$ with $\se^{\prime} < \se^{\prime \prime}$, \be
\la{TwoStopgenb}
 \Pr \li \{ \ba < \Up_n (\bs{\mcal{X}}_n; \se^{\prime}, \se^{\prime \prime}) < \f{1}{\al} \mid \se \ri \}  \to  0
\ee as the sample number $n$ tends to $N^*$.

 (II) For arbitrary $\al \in (0, 1)$ and $\se, \se^{\prime}, \se^{\prime \prime} \in \Se$ with
$\se^{\prime} < \se^{\prime \prime} \leq \se$, \bee \la{TwoStopgenbb} \Pr \li \{ \Up_n (\bs{\mcal{X}}_n; \se^{\prime}, \se^{\prime \prime}) \geq
\f{1}{\al} \mid \se \ri \}  \to  1 \eee as the sample number $n$ tends to $N^*$.

(III) For arbitrary $\ba \in (0, 1)$ and $\se, \se^{\prime}, \se^{\prime \prime} \in \Se$ with $\se \leq \se^{\prime} < \se^{\prime \prime}$,
\bee \la{TwoStopgenbbc}
 \Pr \li \{ \Up_n (\bs{\mcal{X}}_n; \se^{\prime}, \se^{\prime \prime}) \leq \ba \mid \se \ri \}  \to  1
\eee as the sample number $n$ tends to $N^*$.

\eeT

Theorem \ref{genSTgen} can be established by mimicking the argument of the termination property of Theorem \ref{MSPRT_Composite} as in Appendix
\ref{MSPRT_Composite_app}.

It should be noted that (\ref{TwoStopgenb}) implies  \be \la{ineqlater} \Pr \li \{ \ba < \Up_n (\bs{\mcal{X}}_n; \se^{\prime}, \se^{\prime
\prime}) < \f{1}{\al} \; \tx{for all} \; n \in \mcal{N} \mid \se \ri \} = 0, \ee  which has been established in \cite[Appendix A.1]{Wald}, under
a very general assumption,  for the termination property of Wald's sequential probability ratio tests on two hypotheses.  However,
(\ref{ineqlater}) does not imply (\ref{TwoStopgenb}).

 Actually, in the case that $X_1, X_2, \cd$ are i.i.d samples of $X$ parameterized by $\se \in \Se$, the assumption (I) of
 Theorem \ref{genSTgen} holds under fairy general conditions, as can be seen by the following result.

\beT
 \la{genST}

Let $\se, \se^{\prime}, \se^{\prime \prime} \in \Se$.  Assume that $\Pr \{  f (X; \; \se^{\prime} )  f (X; \; \se^{\prime \prime}) = 0 \mid \se
\}  = 0$ and that the variance of $\ln \f{ f (X; \; \se^{\prime \prime}) } { f (X; \; \se^{\prime} ) }$ is positive and finite. Then, for
arbitrary $\al, \ba \in (0, 1)$, \be \la{TwoStopgen}
 \lim_{n \to \iy} \Pr \li \{ \ba < \Up_n (\bs{\mcal{X}}_n; \se^{\prime}, \se^{\prime \prime}) < \f{1}{\al} \mid \se \ri \}  =  0.
\ee
 \eeT

See Appendix \ref{genST_app} for a proof.  It should be noted that if $\Pr \{  f (X; \; \se^{\prime} )   f (X; \; \se^{\prime \prime}) = 0 \mid
\se \} > 0$, then $\lim_{n \to \iy} \Pr \li \{ \ba < \Up_n (\bs{\mcal{X}}_n; \se^{\prime}, \se^{\prime \prime}) < \f{1}{\al} \mid \se \ri \}  =
0$.

In the case that $X_1, X_2, \cd$ are i.i.d samples of $X$ parameterized by $\se \in \Se$, the assumptions (II) and (III) of Theorem
\ref{genSTgen} are valid under fairy general conditions, as can be seen by the following result.

\beT  \la{STChe}  Let $\se, \se^{\prime}, \se^{\prime \prime} \in \Se$ and $\al, \ba \in (0, 1)$.  Assume that $\Pr \{  f (X; \; \se^{\prime} )
f (X; \; \se^{\prime \prime}) = 0 \mid \se \}  = 0$ and that the variances of $\ln f(X; \se^{\prime})$ and $\ln f(X; \se^{\prime \prime})$
associated with $\se$ are positive and finite.  Then,

(I) $ \lim_{n \to \iy} \Pr \li \{ \Up_n (\bs{\mcal{X}}_n; \se^{\prime}, \se^{\prime \prime}) \geq \f{1}{\al} \mid \se \ri \}  =  1$ holds under
the additional assumption that $\bb{E} [ \ln f(X; \se^{\prime}) \mid \se ] < \bb{E} [ \ln f(X; \se^{\prime \prime}) \mid \se ]$.

(II) $ \lim_{n \to \iy} \Pr \li \{ \Up_n (\bs{\mcal{X}}_n; \se^{\prime}, \se^{\prime \prime}) \leq \ba \mid \se \ri \}  =  1$ holds under the
additional assumption that $\bb{E} [ \ln f(X; \se^{\prime}) \mid \se ]
> \bb{E} [ \ln f(X; \se^{\prime \prime}) \mid \se ]$.

\eeT

See Appendix \ref{STChe_app} for a proof.

It should be noted that if $\Pr \{  f (X; \; \se^{\prime} ) = 0 \mid \se \}  > 0$, then $ \lim_{n \to \iy} \Pr \{ \Up_n (\bs{\mcal{X}}_n;
\se^{\prime}, \se^{\prime \prime}) \geq \f{1}{\al} \mid \se \}  = 1$. Similarly, if $\Pr \{  f (X; \; \se^{\prime \prime}) = 0 \mid \se \} > 0$,
then $ \lim_{n \to \iy} \Pr \li \{ \Up_n (\bs{\mcal{X}}_n; \se^{\prime}, \se^{\prime \prime}) \leq \ba \mid \se \ri \}  =  1$.

\subsection{One-sided Hypotheses}

It should be noted that in the special context of testing two hypotheses, our CSPRTs reduce to Wald's SPRTs.

For the problem of testing simple hypotheses $\mscr{H}_0: \se = \se_0$ versus $\mscr{H}_1: \se = \se_1$, the likelihood function $f_n
(\bs{\mcal{X}}_n; \se)$ is unimodal with respect to $\se \in \Se$, since there are only two values in the parameter space $\Se$.  Therefore, the
required assumption of our CSPRT is the same as that of Wald's SPRT.

For the problem of testing composite hypotheses  $\mscr{H}_0: \se \leq \se_0$ versus $\mscr{H}_1: \se \geq \se_1$,  our CSPRT requires the
assumption that the test will surely eventually terminate and that the likelihood function $f_n (\bs{\mcal{X}}_n; \se)$ is unimodal with respect
to $\se \in \Se$.  It has been previously known that the SPRT is applicable to the composite hypotheses under the assumption that the SPRT will
surely eventually terminate and that the relevant likelihood ratio is monotone.

We would like to point out that there are some situations where the relevant likelihood ratio does not possess the monotonicity property, but
the likelihood function $f_n (\bs{\mcal{X}}_n; \se)$ is unimodal with respect to $\se \in \Se$.  To illustrate, consider hypotheses regarding
the distribution of random variable $X$ uniformly distributed on $[\se - a, \se + a]$ with known $a
> 0$ and unknown parameter $\se$.  Suppose one wish to test hypotheses on $\se$ based on i.i.d. samples $X_1, X_2, \cd$ of $X$.
Since for any sample number $n$, the likelihood ratio needs to be expressed in terms of $\min \{ X_1, \cd, X_n \}$ and $\max \{ X_1, \cd, X_n
\}$, we can conclude that the likelihood ratio  does not possess the monotonicity property.  However, it can be readily shown that the
likelihood function $f_n (\bs{\mcal{X}}_n; \se)$ is unimodal with respect to $\se \in \Se$.  From this discussion, it can be seen that our
result in Theorem \ref{MSPRT_Composite} has extended the applications of Wald's SPRTs to a wider variety  of composite hypotheses.

\subsection{Two-sided Hypotheses}

Consider a classical problem of testing two-sided hypotheses $H_0: \se = \vse_0$ versus $H_1: \se \neq \vse_0$, with $\vse_0 \in \Se$.  As
pointed out by Wald \cite[Section 4.4.4, page 77]{Wald}, it is a common contention that the acceptance of $H_0$ will not be considered a serious
error if $\se \neq \vse_0$ but is near $\vse_0$. However, there will be, in general, two parameter values $a$ and $b$ with $a < \vse_0 < b$ such
that the acceptance of $H_0$ is considered an error of practical importance if (and only if) $\se \notin (a, b)$.  Thus, the region of
preference for rejection may be defined as the set of all values $\se$ for which $\se \notin (a, b)$.  The region of preference for acceptance
will consist of the single value $\vse_0$, and the region of indifference will be the set of all values $\se$ for which $(a, \vse_0) \cup
(\vse_0, b)$. To control the risk of making wrong decision, it is typically required that \be \la{twoa}
 \Pr \{ \tx{Reject} \; H_0 \mid \se \} \leq \al \qu
\tx{for} \; \se = \vse_0 \ee and \be \la{threb} \Pr \{ \tx{Accept} \; H_0 \mid \se \} \leq \ba \qu \tx{for} \; \se \in \Se \; \tx{such that} \;
\se \notin (a, b). \ee To solve this problem, Wald proposed the principle of weight function.  However, an appropriate weight function is
difficult to find, especially for discrete distributions.   We propose to solve the problem by constructing CSPRT for the following three new
hypotheses
\[
\mscr{H}_0: \se \leq \f{a + \vse_0}{2}, \qqu \mscr{H}_1: \f{a + \vse_0}{2} < \se \leq \f{b + \vse_0}{2}, \qqu \mscr{H}_2:  \se
> \f{b + \vse_0}{2}
\]
so that \bee  \Pr \{ \tx{Reject} \; \mscr{H}_0 \mid \se \} \leq \f{\ba}{2}  \;  \tx{for} \; \se \leq a; \qu  \Pr \{ \tx{Reject} \; \mscr{H}_1
\mid \se \} \leq \al  \; \tx{for} \; \se = \vse_0; \qu  \Pr \{ \tx{Reject} \; \mscr{H}_2 \mid \se \} \leq \f{\ba}{2}  \; \tx{for} \; \se \geq b.
\eee This can be accomplished by applying the
CSPRT with $m = 3, \; \de_0 = \f{\ba}{2}, \; \de_1 = \al, \; \de_2 = \f{\ba}{2}$ and
\bee &  & \se_1 = \f{a + \vse_0}{2}, \qqu \se_2 = \f{b + \vse_0}{2},\qqu \se_1^\prime = a, \qqu \se_1^{\prime \prime} = \vse_0, \qqu \se_2^\prime = \vse_0, \qqu \se_2^{\prime \prime} = b,\\
&  & \al_1 = \f{\ba}{2}, \qqu \al_2  = \f{\al}{2}, \qqu \al_3 = 0,\qqu \ba_0 = 0, \qqu \ba_1 = \f{\al}{2}, \qqu \ba_2 = \f{\ba}{2}. \eee

At the termination of the CSPRT, the decision on the original hypotheses $H_0$ versus $H_1$ is made based on the decision on the new hypotheses
$\mscr{H}_0, \; \mscr{H}_1$ and $\mscr{H}_2$ by the following rule:

Accept $H_0$ if $\mscr{H}_1$ is accepted; Reject $H_0$ if either $\mscr{H}_0$ or $\mscr{H}_2$ is accepted.

Based on this proposal, it can be readily shown that the risk requirements (\ref{twoa}) and (\ref{threb}) are satisfied.

\subsection{CSPRTs on Parameters of Exponential Family}

In this section, we shall show that the CSPRTs can be applied to the parameters of the exponential family under mild assumptions.   Let $X$ be a
random variable with PDF or PMF of the form
\[
f_X (x; \se) = h (x) \exp [ u (\se) T(x) - v (\se) ],
\]
where $T(x)$ and $h(x)$ are functions of $x$, and $u (\se), \; v (\se)$ are functions of $\se \in \Se$.   We have obtained the following
results.

 \beT \la{expST}
Assume that  $\f{ d v (\se) }{d \se} = \se \f{ d u (\se) }{d \se}$ and that $\f{ d u (\se) }{d \se} > 0$ for $\se \in \Se$.  Let $X_1, X_2, \cd$
be i.i.d. samples of $X$.  Then, for any  $n \in \bb{N}$, the likelihood function $f_n ( \bs{\mcal{X}}_n; \se)$ is unimodal with respect $\se
\in \Se$. Moreover, for arbitrary $\al, \ba \in (0, 1)$ and $\se, \; \se^{\prime}, \; \se^{\prime \prime} \in \Se$ with $\se^{\prime} <
\se^{\prime \prime}$,
\[
\lim_{n \to \iy} \Pr \li \{ \ba < \Up_n ( \bs{\mcal{X}}_n; \se^{\prime}, \se^{\prime \prime} )  < \f{1}{\al} \mid \se  \ri \} = 0.
\]
\eeT

See Appendix \ref{expST_app} for a proof.  It can be readily verified that the assumption of Theorem \ref{expST} is satisfied for the binomial,
Poisson,  normal, exponential, gamma, geometric and negative binomial distributions.

\subsection{CSPRTs on Proportion of Finite Population}

Consider a finite population of $N$ units among which there are $Np$ units having a certain attribute, where $p \in \Se \DEF \{ \f{i}{N}: i = 0,
1, \cd, N \}$.   Many practical problems can be formulated as the multiple hypotheses testing problem defined by (\ref{mainpr}) and
(\ref{mainreq}), with the parameter $\se$ identified as $p$.  For such a problem, consider sampling without replacement. As before, define a
Bernoulli random variable $X_n$ such that $X_n$ assumes values $1$  or $0$ in accordance with whether the $n$-th drawn unit has the attribute.
This leads to a sequence of dependent Bernoulli random variables $X_1, \cd, X_N$ parameterized by $p \in \Se$.  The following analysis shows
that our CSPRTs can be applied to the general multiple hypotheses testing problem.

Clearly, the likelihood function is
\[
f_n ( \bs{\mcal{X}}_n; p ) = \f{ \bi{N p}{ K_n } \bi{N - N p }{n - K_n}  } { \bi{n}{ K_n } \bi{N}{ n }  },
\]
where $K_n = \sum_{i=1}^n X_i$.  Let $\al, \ba \in (0,1)$ and $p, p^{\prime}, p^{\prime \prime} \in \Se$ with $p^{\prime} < p^{\prime \prime}$.

In the case of $p \leq p^{\prime}$, we have $f_N ( \bs{\mcal{X}}_N; p^{\prime \prime}) = 0$. Thus,  $f_N ( \bs{\mcal{X}}_N; p^{\prime \prime})
> \ba f_N ( \bs{\mcal{X}}_N; p^{\prime})$ is violated.

In the case of $p \geq p^{\prime \prime}$, we have $f_N ( \bs{\mcal{X}}_N; p^{\prime}) = 0$. Thus,  $f_N ( \bs{\mcal{X}}_N; p^{\prime})
> \al f_N ( \bs{\mcal{X}}_N; p^{\prime \prime})$ is violated.

In the case of $p^{\prime} < p < p^{\prime \prime}$, it must be true that $f_N ( \bs{\mcal{X}}_N; p^{\prime}) = f_N ( \bs{\mcal{X}}_N; p^{\prime
\prime}) = 0$, which implies that both $f_N ( \bs{\mcal{X}}_N; p^{\prime \prime}) > \ba f_N ( \bs{\mcal{X}}_N; p^{\prime})$ and  $f_N (
\bs{\mcal{X}}_N; p^{\prime})
> \al f_N ( \bs{\mcal{X}}_N; p^{\prime \prime})$ are violated.

This proves that
\[
\Pr \li \{ \ba < \Up_n ( \bs{\mcal{X}}_n; p^{\prime}, p^{\prime \prime} ) < \f{1}{\al} \mid p \ri \} \to 0
\]
as $n \to N$.  It can be shown by direct computation that $f_n ( \bs{\mcal{X}}_n; p )$ is unimodal with respect to $p$.

\subsection{Unimodal Property of Various Distributions}

In addition to the exponential family and the distribution associated with a sampling without replacement from a finite population, the
likelihood functions of a wide variety of distributions have the desired unimodal properties which permit the applications of CSPRTs.  A few of
such distributions are outlined in the sequel.

\subsubsection{Positive Power Law Distribution}

A random variable $X$ is said to have a positive power law distribution if the density function of $X$ is given by
\[
f_X (x; \ga, \ka) = \bec \f{\ka+1}{\ga^{\ka +1}} x^\ka  & \tx{for} \; x \in [0, \ga], \\
0 & \tx{for} \; x \notin [0, \ga], \eec
\]
where $\ka \geq 0$ and $\ga > 0$.  Clearly, taking $\ka = 0$ gives the uniform distribution.   It can be checked that for a given $\ga > 0$, the
likelihood function $f_n( \bs{\mcal{X}}_n; \ga, \ka)$ is unimodal with respect to $\ka$.  On the other hand, when $\ka \geq 0$ is fixed, $f_n(
\bs{\mcal{X}}_n; \ga, \ka)$ is unimodal with respect to $\ga > 0$.

\subsubsection{Pareto Distribution}

The Pareto distribution is given in density-function form by
\[
f_X (x; \ga, \ka) = \bec \f{\ka}{\ga} \li ( \f{\ga}{x} \ri )^{\ka + 1}  & \tx{for} \; x \in [\ga, \iy), \\
0 & \tx{for} \; x \notin [\ga, \iy), \eec
\]
where $\ka > 0$ and $\ga > 0$.  It can be shown that for any given $\ga > 0$, $f_n( \bs{\mcal{X}}_n; \ga, \ka)$ is unimodal with respect to
$\ka$. When $\ka$ is fixed, $f_n( \bs{\mcal{X}}_n; \ga, \ka)$ is unimodal with respect to $\ga > 0$.

\subsubsection{Normal Distribution with Known Mean}

The normal distribution is given in density-function form by
\[
f_X (x; \mu, \si) = \f{1}{\sq{2 \pi} \si} \exp \li ( - \f{(x - \mu)^2}{2 \si^2} \ri ),
\]
where $- \iy < \mu < \iy$ and $\si > 0$. It can be shown that for any given $\mu$, the likelihood function $f_n( \bs{\mcal{X}}_n; \mu, \si)$ is
unimodal with respect to $\si$.

\subsubsection{Laplace Distribution}

A random variable $X$ is said to have a Laplace distribution if the density function of $X$ is given by
\[
f_X (x; \mu, \nu) = \f{1}{2 \nu} \exp \li ( - \f{|x - \mu|}{\nu} \ri ),
\]
where $- \iy < \mu < \iy$ and $\nu > 0$. It can be shown that for any given $\mu$, the likelihood function $f_n( \bs{\mcal{X}}_n; \mu, \nu)$ is
unimodal with respect to $\nu$.

\subsubsection{Negative Exponential Distribution}

The negative exponential distribution is given in density-function form by
\[
f_X(x; \mu, \nu) = \bec \f{1}{\nu} \exp \li (  - \f{x - \mu}{\nu} \ri )  & \tx{for} \; x \in [\mu, \iy), \\
0 & \tx{for} \; x \notin [\mu, \iy), \eec
\]
where $- \iy < \mu < \iy$ and $\nu > 0$.  Clearly, for any given $\mu$, the likelihood function $f_n( \bs{\mcal{X}}_n; \mu, \nu)$ is unimodal
with respect to $\nu > 0$.  On the other hand, when $\nu > 0$ is fixed, $f_n( \bs{\mcal{X}}_n; \mu, \nu)$ is unimodal with respect to $\mu$.

\subsubsection{Weibull Distribution}

The Weibull distribution is given in density-function form by
\[
f_X (x; \lm, \ka) = \f{\ka}{\lm} \li ( \f{x}{\lm} \ri )^{\ka - 1} \exp \li ( - \li ( \f{x}{\lm} \ri )^\ka \ri ), \qqu x > 0, \qu \ka > 0, \qu
\lm > 0
\]
It can be shown that for any given $\ka > 0$, the likelihood function $f_n( \bs{\mcal{X}}_n; \lm, \ka)$ is unimodal with respect to $\lm > 0$.

\section{Continuous-Time Stochastic Processes} \la{secCon}

By a similar approach as that of the CSPRTs for the discrete-time process $(X_n)_{n \in \mcal{N}}$, we can develop CSPRTs for a continuous-time
processes $(X_t)_{t \in [0, \iy)}$ parameterized by $\se \in \Se$. Throughout Sections \ref{sec3A} and \ref{sec8A}, let $(X_t)_{t \in [0, \iy)}$
be a right-continuous stochastic process parameterized by $\se \in \Se$ and let the probability mass or density function of $X_t$ be denoted by
$f_t(. ; \se)$ for $t \in [0, \iy)$.  Assume that $f_t(x; \se)$ is right-continuous with respect to $t \in [0, \iy)$ for any $\se \in \Se$ and
$x \in \bb{R}$.

\subsection{Maximal Inequality and Confidence Sequences}  \la{sec3A}

For parameter values $\se^\prime, \se^{\prime \prime} \in \Se$, define likelihood ratio $\Up_t (X_t; \se^\prime, \se^{\prime \prime}) = \f{ f_t
(X_t; \se^{\prime \prime}) } { f_t (X_t; \se^\prime)  }$ for $t \in [0, \iy)$.  We have established the following results on maximal
inequalities and confidence sequences.

\beT  \la{CI sequence Continuous} Assume that for arbitrary integer $n$ and real numbers $t_i, \; i = 0, \cd, n$ with $0 = t_0 < t_1 < \cd <
t_{n-1} < t_n = t$, the conditional probability mass or density function of $X_{t_i}, \; i = 0, 1, \cd, n - 1$ given the value of $X_t$ does not
depend on $\se$.  Let $\se_0, \se_1 \in \Se$ and $\de \in (0, 1)$.  Then, \be \la{maxcon} \Pr \li \{ \Up_t (X_t; \se_0, \se_1) >  \f{1}{\de} \;
\tx{for some} \; t \in [0, \iy) \mid \se_0 \ri \} \leq \de. \ee Moreover, $\Pr \{ L_t (X_t) \leq \se \; \tx{for all} \; t \mid \se \} \geq 1 -
\f{\de}{2}, \; \;
 \Pr \{  U_t (X_t) \geq \se \; \tx{for all} \; t  \mid \se \} \geq 1 - \f{\de}{2}$
 and $\Pr \{ L_t (X_t) \leq \se \leq U_t (X_t) \; \tx{for all} \; t  \mid \se \} \geq 1 - \de$ for all $\se \in \Se$,
 where $L_t (X_t) = \inf  \{ \vse \in \Se: \Up_t (X_t; \se_1, \vse) \geq \f{\de}{2} \}$ and $U_t (X_t) = \sup \{  \vse \in \Se:
 \Up_t (X_t; \se_0, \vse) \geq \f{\de}{2} \}$.

\eeT

See Appendix \ref{CI sequence Continuous_app} for a proof.  If the likelihood function $f_t (X_t; \se)$ is unimodal with respect to $\se \in
\Se$, then there exists an estimator $\wh{\se}_t$ for $\se$ such that $f_t (X_t; \se)$ is non-decreasing with respect to $\se \in \Se$ no
greater than $\wh{\se}_t$ and is non-increasing with respect to $\se \in \Se$ no less than $\wh{\se}_t$.  Hence, it must be true that
$\{\wh{\se}_t \leq \se_0 \} \subseteq \{ \Up_t (X_t; \se_1, \se_0) \geq 1 \}$ and consequently, $\{ \Up_t (X_t; \se_1, \se_0) < \f{\de}{2}, \;
\wh{\se}_t \leq \se_0 \} \subseteq \{ \Up_t (X_t; \se_1, \se_0) < 1, \; \wh{\se}_t \leq \se_0 \} = \emptyset$. It follows that \bel \li \{ \Up_t
(X_t; \se_1, \se_0) < \f{\de}{2} \ri \} & = & \li \{ \Up_t (X_t; \se_1,
\se_0) < \f{\de}{2}, \; \wh{\se}_t  < \se_0 \ri \} \nonumber\\
& \subseteq & \li \{ \Up_t (X_t; \se_1,  \se) < \f{\de}{2} \; \tx{for all} \;  \se \leq \se_0 \ri \} \la{useulecc} \\
& \subseteq & \{ L_t (X_t) \geq \se_0 \},  \nonumber \eel where (\ref{useulecc}) is also a consequence of the assumption that the likelihood
function $f_t (X_t; \se)$ is unimodal with respect to $\se \in \Se$.

\subsection{CSPRTs on Multiple Hypotheses} \la{sec8A}

For the multi-hypotheses testing problem defined by (\ref{mainpr}) and (\ref{mainreq}), we propose a CSPRT with stopping and decision rules as
follows:

\vspace{0.05in}

\begin{tabular} {|l |}
\hline $ \tx{ {\it Continue observing $(X_t)_{t \in [0, \iy)}$ until there exists an index $j$ in the set $\{0, 1, \cd, m - 1 \}$ such that}}$\\
$\tx{ {\it $\Up_t (X_t; \se_i^{\prime}, \se_i^{\prime \prime}) > \f{1}{\al_i}$ for $0 < i \leq j$ and $\Up_t (X_t; \se_i^{\prime}, \se_i^{\prime
\prime}) < \ba_i$ for $j < i < m$.
At the termination}}$\\
$\tx{ {\it of the observational procedure, accept $\mscr{H}_j$ with the index $j$ satisfying the stopping condition}}$.
\\ \hline
\end{tabular}

\vspace{0.05in}

We have established that the above CSPRT has the following properties.

\beT \la{MSPRT_Composite_Continuous}
 Assume that for arbitrary integer $n$ and real numbers $t_i, \; i = 0, \cd, n$ with $0 = t_0 < t_1 < \cd <
t_{n-1} < t_n = t$, the conditional probability mass or density function of $X_{t_i}, \; i = 0, 1, \cd, n - 1$ given the value of $X_t$ does not
depend on $\se$.  Assume that the likelihood function $f_t (X_t; \se)$ is unimodal with respect to $\se \in \Se$ for any positive number $t$.
If the observational process will eventually terminate according to the stopping rule with probability $1$, then the following statements
(I)--(III) hold true:

(I) $\Pr \{ \tx{Reject} \; \mscr{H}_i \mid \se \} \leq \al_{i + 1} + \ba_i$ for $0 \leq i \leq m - 1$ and $\se \in \varTheta_i$.

(II) For $j = 1, \cd, m - 1$, $\Pr \{ \tx{Accept} \; \mscr{H}_i \; \tx{with some index $i$ no less than $j$} \mid \se \}$ is no greater than
$\al_j$ and is non-decreasing with respect to $\se \in \Se$ no greater than $\se_j^\prime$.

(III) For $j = 1, \cd, m - 1$, $\Pr \{ \tx{Accept} \; \mscr{H}_i \; \tx{with some index $i$ less than $j$} \mid \se \}$ is no greater than
$\ba_j$ and is non-increasing with respect to $\se \in \Se$ no less than $\se_j^{\prime \prime}$.

Moreover, the sampling process will eventually terminate according to the stopping rule with  probability $1$,  provided that the following
additional assumption is satisfied: For arbitrary $\al, \ba \in (0, 1)$ and $\se \in \Se$,
\[
\lim_{t \to \iy} \Pr \li \{ \ba \leq \Up_t (X_t; \se_i^{\prime}, \se_i^{\prime \prime}) \leq  \f{1}{\al}  \mid \se \ri \} = 0, \qqu i = 1, \cd,
m - 1
\] \eeT

The proof of Theorem \ref{MSPRT_Composite_Continuous}  is similar to that of Theorem \ref{MSPRT_Composite}.

For testing simple hypothesis defined by (\ref{simdefA}) and (\ref{simdefB}), we propose a CSPRT with stopping and decision rules as follows:

\vspace{0.05in}

\begin{tabular} {|l |}
\hline $ \tx{ {\it Continue observing $(X_t)_{t \in [0, \iy)}$ until there exists an index $j$ in the set $\{0, 1, \cd, m - 1 \}$ such that}}$\\
$\tx{ {\it $\Up_t (X_t; \se_{i-1}, \se_{i}) > \f{1}{\al_i}$ for $0 < i \leq j$ and $\Up_t (X_t; \se_{i-1}, \se_{i}) < \ba_i$ for $j < i <
m$. At the termination}}$\\
$\tx{ {\it of the observational procedure, accept $\mscr{H}_j$ with the index $j$ satisfying the stopping condition}}$.
\\ \hline
\end{tabular}

\vspace{0.05in}

We have established that such CSPRT possesses the following properties.

\beT  \la{MSPRT_Simple_Continuous}  Assume that for arbitrary integer $n$ and real numbers $t_i, \; i = 0, \cd, n$ with $0 = t_0 < t_1 < \cd <
t_{n-1} < t_n = t$, the conditional probability mass or density function of $X_{t_i}, \; i = 0, 1, \cd, n - 1$ given the value of $X_t$ does not
depend on $\se$.  If the observational process will eventually terminate according to the stopping rule with probability $1$, then $\Pr \{
\tx{Reject} \; \mscr{H}_i \mid \se_i \} \leq \al_{i+1} + \ba_i$ for $0 \leq i \leq m - 1$.  Moreover, the observational process will eventually
terminate according to the stopping rule with probability $1$,  provided that the likelihood function $f_t (X_t; \se_i)$ is unimodal with
respect to $\se \in \Se$,  and that for arbitrary $\al, \ba \in (0, 1)$ and $\se \in \Se$,
\[
\lim_{t \to 0} \Pr \li \{ \ba \leq \Up_t (X_t; \se_{i-1}, \se_{i}) \leq  \f{1}{\al} \mid \se \ri \} = 0, \qqu i = 1, \cd, m - 1. \]
 \eeT

The proof of Theorem \ref{MSPRT_Simple_Continuous}  is similar to that of Theorem \ref{MSPRT_Simple}.

\subsection{CSPRTS on Arrival Rates of Poisson Processes}

Consider a Poisson process $(X_t)_{t \in [0, \iy)}$ with an arrival rate $\lm > 0$.  Note that for $\ga > 0$,
\[
\li \{ \f{ f_t (X_t; \lm_1) } { f_t (X_t; \lm_0) } > \ga \ri \}  = \li \{ X_t > \f{ (\lm_1 - \lm_0) t + \ln \ga  } {\ln \f{\lm_1}{\lm_0} } \ri
\}.
\]
For testing multiple composite hypotheses defined by (\ref{mainpr}) and (\ref{mainreq}), with $\se, \; \se_i^\prime, \; \se_i, \; \se_i^{\prime
\prime}$
 identified as $\lm, \; \lm_i^\prime, \; \lm_i, \; \lm_i^{\prime \prime}$ respectively, we propose a CSPRT with stopping and decision rules as
 follows:

\vspace{0.05in}

\begin{tabular} {|l |}
\hline $ \tx{ {\it Continue observing $(X_t)_{t \in [0, \iy)}$ until there exists an index $j$ in the set $\{0, 1, \cd, m - 1 \}$ such that}}$\\
$\tx{ {\it $X_t > \f{ (\lm_i^{\prime \prime} - \lm_i^{\prime}) t + \ln \f{1}{\al_i}  } {\ln \f{\lm_i^{\prime \prime}}{\lm_i^{\prime}} }$ for $0
< i \leq j$ and $X_t < \f{ (\lm_i^{\prime \prime} - \lm_i^{\prime}) t + \ln \ba_i  } {\ln \f{\lm_i^{\prime \prime}}{\lm_i^{\prime}} }$ for $j <
i < m$.  At the termination of}}$\\
$\tx{ {\it the observational procedure, accept $\mscr{H}_j$ with the index $j$ satisfying the stopping condition.}}$
\\ \hline
\end{tabular}

\vspace{0.05in}

Regarding the above CSPRT, we have shown the following result.

\beT

\la{PoST}

The observational process will eventually terminate according to the stopping rule with probability $1$.  Moreover, the following statements
(I)--(III) hold true:

(I) $\Pr \{ \tx{Reject} \; \mscr{H}_i \mid \lm \} \leq \al_{i + 1} + \ba_i$ for $0 \leq i \leq m - 1$ and $\lm \in \varTheta_i$.

(II) For $j = 1, \cd, m - 1$, $\Pr \{ \tx{Accept} \; \mscr{H}_i \; \tx{with some index $i$ no less than $j$} \mid \lm \}$ is no greater than
$\al_j$ and is non-decreasing with respect to $\lm \in \Se$ no greater than $\lm_j^\prime$.

(III) For $j = 1, \cd, m - 1$, $\Pr \{ \tx{Accept} \; \mscr{H}_i \; \tx{with some index $i$ less than $j$} \mid \lm \}$ is no greater than
$\ba_j$ and is non-increasing with respect to $\lm \in \Se$ no less than $\lm_j^{\prime \prime}$.

\eeT

See Appendix \ref{PoST_app} for a proof.

 For testing multiple simple hypotheses defined by (\ref{simdefA}) and (\ref{simdefB}), with $\se, \; \se_i$
 identified as $\lm, \; \lm_i$ respectively, we propose a CSPRT with stopping and decision rules as
 follows:

\vspace{0.05in}

\begin{tabular} {|l |}
\hline $ \tx{ {\it Continue observing $(X_t)_{t \in [0, \iy)}$ until there exists an index $j$ in the set $\{0, 1, \cd, m - 1 \}$ such that}}$\\
$\tx{ {\it $X_t > \f{ (\lm_i - \lm_{i-1}) t + \ln \f{1}{\al_i}  } {\ln \f{\lm_i}{\lm_{i-1}} }$ for $0 < i \leq j$ and $X_t < \f{ (\lm_i -
\lm_{i-1}) t + \ln \ba_i  } {\ln \f{\lm_i}{\lm_{i-1}} }$ for $j <
i < m$.  At the termination}}$\\
$\tx{ {\it of the observational procedure, accept $\mscr{H}_j$ with the index $j$ satisfying the stopping condition.}}$
\\ \hline
\end{tabular}

\vspace{0.05in}

Regarding the above CSPRT, we have shown the following result.

\beT  \la{MSPRT_Simple_Continuous_Pos}  The observational process will eventually terminate according to the stopping rule with probability $1$.
Moreover,  $\Pr \{ \tx{Reject} \; \mscr{H}_i \mid \lm_i \} \leq \al_{i+1} + \ba_i$ for $0 \leq i \leq m - 1$.
 \eeT

Theorem \ref{MSPRT_Simple_Continuous_Pos}  is a direct consequence of Theorem \ref{MSPRT_Simple_Continuous}.

\subsection{CSPRTS on Parameters of Brownian Motions}

Consider a Brownian  motion $(X_t)_{t \in [0, \iy)}$ with unknown drift $\mu$ and known variance $\si^2$ per unit time. Note that for $\ga > 0$,
\[
\li \{ \f{ f_t (X_t; \mu_1, \si) } { f_t (X_t; \mu_0, \si) } > \ga \ri \} =  \li \{ X_t > \f{(\mu_0 + \mu_1) t}{2} + \f{\si^2}{\mu_1 - \mu_0}
\ln \ga \ri \}.
\]
For testing multiple composite hypotheses defined by (\ref{mainpr}) and (\ref{mainreq}) with $\se, \; \se_i^\prime, \; \se_i, \; \se_i^{\prime
\prime}$
 identified as $\mu, \; \mu_i^\prime, \; \mu_i, \; \mu_i^{\prime \prime}$ respectively, we propose a CSPRT with stopping and decision rules as
 follows:

\vspace{0.05in}

\begin{tabular} {|l |}
\hline $ \tx{ {\it Continue observing $(X_t)_{t \in [0, \iy)}$ until there exists an index $j$ in the set $\{0, 1, \cd, m - 1 \}$ such}}$\\
$\tx{ {\it that $X_t > \f{(\mu_i^{\prime} + \mu_i^{\prime \prime}) t}{2} + \f{\si^2}{\mu_i^{\prime \prime} - \mu_i^{\prime}} \ln \f{1}{\al_i}$
for $0 < i \leq j$ and $X_t < \f{(\mu_i^{\prime} + \mu_i^{\prime \prime}) t}{2} + \f{\si^2}{\mu_i^{\prime \prime} - \mu_i^{\prime}} \ln \ba_i$
for $j < i < m$.}}$\\
$\tx{ {\it At the termination of the observational procedure, accept $\mscr{H}_j$ with the index $j$ satisfying the}}$\\
$\tx{ {\it stopping condition}}$.
\\ \hline
\end{tabular}

\vspace{0.05in}

With regard to above CSPRT,  we have shown the following results.

\beT \la{norST} The observational process will eventually terminate according to the stopping rule with probability $1$.  Moreover, the
following statements (I)--(III) hold true:

(I) $\Pr \{ \tx{Reject} \; \mscr{H}_i \mid \mu \} \leq \al_{i + 1} + \ba_i$ for $0 \leq i \leq m - 1$ and $\mu \in \varTheta_i$.

(II) For $j = 1, \cd, m - 1$, $\Pr \{ \tx{Accept} \; \mscr{H}_i \; \tx{with some index $i$ no less than $j$} \mid \mu \}$ is no greater than
$\al_j$ and is non-decreasing with respect to $\mu \in \Se$ no greater than $\mu_j^\prime$.

(III) For $j = 1, \cd, m - 1$, $\Pr \{ \tx{Accept} \; \mscr{H}_i \; \tx{with some index $i$ less than $j$} \mid \mu \}$ is no greater than
$\ba_j$ and is non-increasing with respect to $\mu \in \Se$ no less than $\mu_j^{\prime \prime}$.  \eeT

See Appendix \ref{norST_app} for a proof.

For testing multiple simple hypotheses defined by (\ref{simdefA}) and (\ref{simdefB}) with $\se, \; \se_i$ identified as $\mu, \; \mu_i$
respectively, we propose a CSPRT with stopping and decision rules as follows:

\vspace{0.05in}

\begin{tabular} {|l |}
\hline $ \tx{ {\it Continue observing $(X_t)_{t \in [0, \iy)}$ until there exists an index $j$ in the set $\{0, 1, \cd, m - 1 \}$ such that}}$\\
$\tx{ {\it $X_t > \f{(\mu_{i-1} + \mu_i) t}{2} + \f{\si^2}{\mu_i - \mu_{i-1}} \ln \f{1}{\al_i}$ for $0 < i \leq j$ and $X_t < \f{(\mu_{i-1} +
\mu_i) t}{2} + \f{\si^2}{\mu_i - \mu_{i-1}} \ln \ba_i$ for $j <
i < m$.}}$\\
$\tx{ {\it At the termination of the observational procedure, accept $\mscr{H}_j$ with the index $j$ satisfying}}$\\
$\tx{ {\it the stopping condition}}$.
\\ \hline
\end{tabular}

\vspace{0.05in}

Same results as in Theorem \ref{MSPRT_Simple_Continuous_Pos} hold for above CSPRT.

\section{Conclusion} \la{secCons}

In this paper, we have established consecutive sequential probability ratio tests for testing multiple statistical hypotheses.  Our tests are
derived based on the principle of probabilistic comparison. Simple analytic formulae are derived for controlling the risk of making wrong
decisions.  We have demonstrated that the new tests can be applied to a wide variety of statistical hypotheses.

\appendix

\sect{Proof of Theorem  \ref{Multi-Valued Inclusion Principle}} \la{Multi-Valued Inclusion Principle_app}

By the assumption that $\cup_{i = 0}^{m-1} \{ \bs{\mcal{L}} = a_i, \; \bs{\mcal{U}} = b_{i + 1} \} = \Om$, we have $\cup_{i = 0}^{m-1} \{
\bs{\mcal{L}} = a_i \} = \Om$ and $\cup_{i = 0}^{m-1} \{ \bs{\mcal{U}} = b_{i + 1} \} = \Om$. Therefore, for $\se \in \varTheta_i$, we have \bee
&  & \{ \bs{\mcal{L}} \geq \se \} = \bigcup_{\ell \in I_{\bs{\tau}} } \{ \bs{\tau} = \ell, \; \bs{\mcal{L}} \geq \se \} = \bigcup_{\ell \in
I_{\bs{\tau}} } \{ \bs{\tau} = \ell, \;
\bs{\mcal{L}} \geq a_{i+1} \} = \{ \bs{\mcal{L}} \geq a_{i+1} \},\\
&  & \{ \bs{\mcal{U}} \leq \se \} = \bigcup_{\ell \in I_{\bs{\tau}} } \{ \bs{\tau} = \ell, \; \bs{\mcal{U}} \leq \se \} = \bigcup_{\ell \in
I_{\bs{\tau}} } \{ \bs{\tau} = \ell, \; \bs{\mcal{U}} \leq b_i \} = \{ \bs{\mcal{U}} \leq b_i \}. \eee

For $i = m - 1$, we have $\varTheta_i = \varTheta_{m-1} = [b_{m - 1}, \iy ), \; L_{\ell, i + 1} = L_{\ell, m} = - \iy, \; a_{i+1} = a_m = \iy$
and hence, $\Pr \{ \bs{\mcal{L}} \geq \se \} = \Pr \{ \bs{\mcal{L}} \geq a_{i+1} \} = 0 = \Pr \{ L_{\ell, i + 1} \geq a_{i+1} \; \tx{for some}
\; \ell \in I_{\bs{\tau}} \} = 0$. As a consequence of the assumption that $\cup_{i = 0}^{m-1} \{ \bs{\mcal{L}} = a_i, \; \bs{\mcal{U}} = b_{i +
1} \} = \Om$, we have \bee &   & \{ \bs{\tau} = \ell, \; \bs{\mcal{L}} = a_j \} = \{ \bs{\mcal{L}} = a_j \} \cap \{ \bs{\tau} = \ell \} = \{
\bs{\mcal{L}} = a_j \} \cap  ( \cup_{i = 0}^{m-1} \{ \bs{\tau} = \ell, \; \bs{\mcal{L}} = a_i, \;
\bs{\mcal{U}} = b_{i + 1} \} )\\
&  & = \{ \bs{\tau} = \ell, \; \bs{\mcal{L}} = a_j, \; \bs{\mcal{U}} = b_{j + 1} \} \eee for $j = 0, 1, \cd, m - 1$.  Hence, for $i = 0, 1, \cd,
m - 2$ and $\se \in \varTheta_i$, we have \bee &  & \{ \bs{\mcal{L}} \geq \se \} = \{ \bs{\mcal{L}} \geq a_{i + 1} \} = \bigcup_{\ell \in
I_{\bs{\tau}} }  \{ \bs{\tau} = \ell, \; \bs{\mcal{L}} \geq a_{i + 1} \} = \bigcup_{\ell \in I_{\bs{\tau}} } \bigcup_{j
> i} \{ \bs{\tau} = \ell, \;
\bs{\mcal{L}} = a_j \}\\
&  & = \bigcup_{\ell \in I_{\bs{\tau}} } \bigcup_{j > i} \{ \bs{\tau} = \ell, \; \bs{\mcal{L}} = a_j, \; \bs{\mcal{U}} = b_{j
+ 1} \}\\
&  & \subseteq  \bigcup_{\ell \in I_{\bs{\tau}} } \bigcup_{j > i} \{ L_{\ell, k} \geq a_k, \; 0 < k \leq j \; \tx{and} \; U_{\ell, k}
\leq b_k, \; j < k < m \}\\
&  &  \subseteq  \bigcup_{\ell \in I_{\bs{\tau}} } \{ L_{\ell, i + 1} \geq a_{i + 1} \} = \{ L_{\ell, i + 1} \geq a_{i+1} \; \tx{for some} \;
\ell \in I_{\bs{\tau}} \}. \eee  For $i = 0$, we have $\varTheta_i = \varTheta_0 = (- \iy, a_1], \; U_{\ell, i} = U_{\ell, 0} = \iy, \; b_i =
b_0 = - \iy$ and hence, $\Pr \{ \bs{\mcal{U}} \leq \se \} = \Pr \{ \bs{\mcal{U}} \leq b_i \} = 0 = \Pr \{ U_{\ell, i} \leq b_i \; \tx{for some}
\; \ell \in I_{\bs{\tau}} \} = 0$.  As a consequence of the assumption that $\cup_{i = 0}^{m-1} \{ \bs{\mcal{L}} = a_i, \; \bs{\mcal{U}} = b_{i
+ 1} \} = \Om$,  we have \bee &   & \{ \bs{\tau} = \ell, \; \bs{\mcal{U}} = b_{j+1} \} = \{
\bs{\mcal{U}} = b_{j+1} \} \cap \{ \bs{\tau} = \ell \} \\
&  & = \{ \bs{\mcal{U}} = b_{j+1} \} \cap  ( \cup_{i = 0}^{m-1} \{ \bs{\tau} = \ell, \; \bs{\mcal{L}} = a_i, \; \bs{\mcal{U}} = b_{i + 1} \} ) =
\{ \bs{\tau} = \ell, \; \bs{\mcal{L}} = a_{j}, \; \bs{\mcal{U}} = b_{j+1} \} \eee for $j = 0, 1, \cd, m - 1$.  Hence, for $i = 1, \cd, m - 1$
and $\se \in \varTheta_i$, we have \bee &  & \{ \bs{\mcal{U}} \leq \se \} = \{ \bs{\mcal{U}} \leq b_i \} = \bigcup_{\ell \in I_{\bs{\tau}} }  \{
\bs{\tau} = \ell, \; \bs{\mcal{U}} \leq b_i \} = \bigcup_{\ell \in I_{\bs{\tau}} } \bigcup_{j < i} \{ \bs{\tau} = \ell, \;
\bs{\mcal{U}} = b_{j+1} \}\\
&  & = \bigcup_{\ell \in I_{\bs{\tau}} } \bigcup_{j < i} \{ \bs{\tau} = \ell, \;
\bs{\mcal{L}} = a_{j}, \; \bs{\mcal{U}} = b_{j + 1} \}\\
&  & \subseteq  \bigcup_{\ell \in I_{\bs{\tau}} } \bigcup_{j < i} \{ L_{\ell, k} \geq a_k, \; 0 < k \leq j \; \tx{and} \; U_{\ell, k}
\leq b_k, \; j < k < m \}\\
&  & \subseteq  \bigcup_{\ell \in I_{\bs{\tau}} } \{ U_{\ell, i} \leq b_{i} \} = \{ U_{\ell, i} \leq b_{i} \; \tx{for some} \; \ell \in
I_{\bs{\tau}} \}. \eee This completes the proof of the theorem.

\section{Proof of Theorem \ref{CI sequence}} \la{CI sequence_app}

We need a preliminary result stated as follows.

\beL \la{Villeinq} Let $\al \in (0, 1)$ and let $\se^\prime, \; \se^{\prime \prime}$ be two parameter values in $\Se$.  Then,
\[
\Pr \li \{ \Up_n (\bs{\mcal{X}}_n; \se^\prime, \se^{\prime \prime} ) \geq \f{1}{\al} \; \tx{for some} \; n \in \mcal{N} \mid \se^\prime \ri \}
\leq \al.
\] \eeL

Actually, the result of Lemma \ref{Villeinq} is due to Ville \cite{Ville}, which was rediscovered by Wald \cite[page 146]{Wald}.

We are now in a position to prove the theorem.  By the definition of the lower confidence limit, we have $\{ L_n (\bs{\mcal{X}}_n) \leq \se_0 \}
\supseteq \li \{ \Up_n (\bs{\mcal{X}}_n; \se_1, \se_0) > \f{\de}{2} \ri \}$.  This implies that  $\{ L_n (\bs{\mcal{X}}_n) > \se_0 \} \subseteq
\li \{ \Up_n (\bs{\mcal{X}}_n; \se_1, \se_0) \leq \f{\de}{2} \ri \}$ and consequently,  $\Pr \{ L_n (\bs{\mcal{X}}_n) > \se \; \tx{for some} \;
n \in \mcal{N} \mid \se \}  \leq \Pr \li \{  \Up_n (\bs{\mcal{X}}_n; \se_1, \se) \leq \f{\de}{2}   \; \tx{for some} \; n \in \mcal{N} \mid \se
\ri \}$ for $\se \in \Se$. It follows from Lemma \ref{Villeinq} that $\Pr \{ L_n (\bs{\mcal{X}}_n) > \se \; \tx{for some} \; n  \in \mcal{N}
\mid \se \} \leq \f{\de}{2}$ for $\se \in \Se$.

Similarly, it follows from the definition of the upper confidence limit that {\small $\{ \Up_n (\bs{\mcal{X}}_n; \se_0,  \se_1) > \f{\de}{2} \}
\subseteq \{ U_n (\bs{\mcal{X}}_n) \geq \se_1 \}$}. This implies that  $\{ U_n (\bs{\mcal{X}}_n) < \se_1 \} \subseteq \{ \Up_n (\bs{\mcal{X}}_n;
\se_0, \se_1) \leq \f{\de}{2} \}$ and consequently, $\Pr \{ U_n (\bs{\mcal{X}}_n) < \se \; \tx{for some} \; n \in \mcal{N} \mid \se \} \leq  \Pr
\{ \Up_n (\bs{\mcal{X}}_n; \se_0, \se) \leq \f{\de}{2}   \; \tx{for some} \; n \in \mcal{N} \mid \se \}$ for $\se \in \Se$.   It follows from
Lemma \ref{Villeinq} that $\Pr \{ U_n (\bs{\mcal{X}}_n) < \se \; \tx{for some} \; n \in \mcal{N} \mid \se \} \leq \f{\de}{2}$ for $\se \in \Se$.
So, by virtue of Bonferroni's inequality, we have $\Pr \{ L_n (\bs{\mcal{X}}_n) \leq \se \leq U_n (\bs{\mcal{X}}_n) \; \tx{for all} \; n  \in
\mcal{N} \mid \se \} \geq 1 - \de$.  This completes the proof of statements (I) and (II).

By the assumption that $f_n ( \bs{\mcal{X}}_n; \se)$ is unimodal with respect to $\se \in \Se$, there exists an estimator $\wh{\se}_n$ of $\se$
such that $f_n ( \bs{\mcal{X}}_n; \se)$ is non-decreasing with respect to $\se \in \Se$ no greater than $\wh{\se}_n$ and is non-increasing with
respect to $\se \in \Se$ no less than $\wh{\se}_n$. Such estimator is referred to as a unimodal-likelihood estimator (ULE) of $\se$.   To show
statement (III), note that as a consequence of the existence of a ULE $\wh{\se}_n$ for $\se$, it must be true that $\{\wh{\se}_n \leq \se_0 \}
\subseteq \{ \Up_n (\bs{\mcal{X}}_n; \se_1, \se_0) \geq 1 \}$ and consequently, $\{ \Up_n (\bs{\mcal{X}}_n; \se_1, \se_0) \leq \f{\de}{2}, \;
\wh{\se}_n \leq \se_0 \} \subseteq \{ \Up_n (\bs{\mcal{X}}_n; \se_1, \se_0) < 1, \; \wh{\se}_n \leq \se_0 \} = \emptyset$. It follows that
{\small \bee \li \{ \Up_n (\bs{\mcal{X}}_n; \se_1, \se_0) \leq \f{\de}{2} \ri \}  =  \li \{ \Up_n (\bs{\mcal{X}}_n; \se_1, \se_0) \leq
\f{\de}{2}, \; \wh{\se}_n  > \se_0 \ri \} \subseteq
 \li \{ \Up_n (\bs{\mcal{X}}_n; \se_1,  \se) \leq \f{\de}{2} \; \tx{for all} \;  \se \leq \se_0 \ri \}
  \subseteq  \{ L_n (\bs{\mcal{X}}_n) \geq \se_0 \}. \eee}
  Similarly, note that as a consequence of the assumption that there exists a ULE $\wh{\se}_n$ for $\se$, it must be true that $\{\wh{\se}_n  \geq
\se_1 \} \subseteq \{ \Up_n (\bs{\mcal{X}}_n; \se_0, \se_1) \geq 1 \}$ and consequently, $\{ \Up_n (\bs{\mcal{X}}_n; \se_0, \se_1) \leq
\f{\de}{2}, \; \wh{\se}_n \geq \se_1 \} \subseteq \{ \Up_n (\bs{\mcal{X}}_n; \se_0, \se_1) < 1, \; \wh{\se}_n \geq \se_1 \} = \emptyset$. It
follows that {\small \bee \li \{ \Up_n (\bs{\mcal{X}}_n; \se_0, \se_1) \leq \f{\de}{2} \ri \}  =  \li \{ \Up_n (\bs{\mcal{X}}_n; \se_0, \se_1)
\leq \f{\de}{2}, \; \wh{\se}_n  < \se_1 \ri \} \subseteq
 \li \{ \Up_n (\bs{\mcal{X}}_n; \se_0,  \se) \leq \f{\de}{2} \; \tx{for all} \;  \se \geq \se_1 \ri \}
  \subseteq  \{ U_n (\bs{\mcal{X}}_n) \leq \se_1 \}. \eee} This completes the proof of the theorem.

\sect{Proof of Theorem \ref{MSPRT_Composite}}  \la{MSPRT_Composite_app}

We need to develop some preliminary results based on the assumptions of the theorem.

 \beL \la{bridge} Let $\al \in (0, 1)$ and let $\se^{\prime} <
\se^{\prime \prime}$ be two parameter values in $\Se$.  Then, \be \la{lemA8} \li \{ \Up_n (\bs{\mcal{X}}_n; \se^{\prime}, \se^{\prime \prime})
\geq \f{1}{\al} \ri \} \subseteq \li \{ \Up_n (\bs{\mcal{X}}_n; \se, \se^{\prime \prime} ) \geq \f{1}{\al} \ri \} \qu \tx{for} \; \; \se \in (-
\iy, \; \se^\prime ] \cap \Se. \ee
 Similarly,
 \be
\la{lemB8}  \{ \Up_n (\bs{\mcal{X}}_n; \se^{\prime}, \se^{\prime \prime}) \leq \al \} \subseteq \{ \Up_n (\bs{\mcal{X}}_n; \se^{\prime}, \se)
\leq \al \} \qu \tx{for} \; \; \se \in [\se^{\prime \prime}, \iy) \cap \Se. \ee \eeL

\bpf

As pointed out in the proof of Theorem \ref{CI sequence} in Appendix \ref{CI sequence_app}, by the assumption that $f_n ( \bs{\mcal{X}}_n; \se)$
is unimodal with respect to $\se \in \Se$, there exists a ULE  $\wh{\se}_n$ for $\se$.

To show (\ref{lemA8}), note that $\{ \Up_n (\bs{\mcal{X}}_n; \se^{\prime}, \se^{\prime \prime}) \geq \f{1}{\al}, \; \wh{\se}_n \leq \se^{\prime}
\} = \emptyset$ and that $\{ \Up_n (\bs{\mcal{X}}_n; \se^\prime, \se^{\prime \prime}) \geq \f{1}{\al}, \; \wh{\se}_n
> \se^{\prime} \} \subseteq
 \{ \Up_n (\bs{\mcal{X}}_n; \se, \se^{\prime \prime}) \geq \f{1}{\al} \}$ for $\se \in   (- \iy, \; \se^\prime ] \cap \Se$.  It follows that
{\small $\{ \Up_n (\bs{\mcal{X}}_n;  \se^{\prime}, \se^{\prime \prime}) \geq \f{1}{\al} \} = \{ \Up_n (\bs{\mcal{X}}_n; \se^{\prime},
\se^{\prime \prime}) \geq \f{1}{\al}, \; \wh{\se}_n >  \se^{\prime} \} \subseteq  \{ \Up_n (\bs{\mcal{X}}_n; \se, \se^{\prime \prime} ) \geq
\f{1}{\al} \}$} for $\se \in (- \iy, \; \se^\prime ] \cap \Se$.

To show (\ref{lemB8}), note that $\{ \Up_n (\bs{\mcal{X}}_n;  \se^{\prime}, \se^{\prime \prime}) \leq \al, \; \wh{\se}_n \geq \se^{\prime
\prime} \} = \emptyset$ and that $\{ \Up_n (\bs{\mcal{X}}_n; \se^{\prime}, \se^{\prime \prime}) \leq \al, \; \wh{\se}_n < \se^{\prime \prime} \}
\subseteq \{ \Up_n (\bs{\mcal{X}}_n; \se^{\prime}, \se) \leq \al \}$ for $\se \in   [\se^{\prime \prime}, \; \iy) \cap \Se$.  It follows that
{\small $\{ \Up_n (\bs{\mcal{X}}_n; \se^{\prime}, \se^{\prime \prime}) \leq \al \}  =  \{ \Up_n (\bs{\mcal{X}}_n; \se^{\prime}, \se^{\prime
\prime}) \leq \al, \; \wh{\se}_n < \se^{\prime \prime} \} \subseteq  \{ \Up_n (\bs{\mcal{X}}_n; \se^{\prime}, \se) \leq \al \}$} for $\se \in
[\se^{\prime \prime}, \; \iy) \cap \Se$.

\epf

\beL \la{lemal}

$\{ \tx{$\mscr{H}_\ell$ with some $\ell > j$ is accepted} \}  \subseteq \{ \Up_n (\bs{\mcal{X}}_n;  \se, \se_{j+1}^{\prime \prime} ) \geq
\f{1}{\al_{j+1}} \; \tx{for some}  \; n \in \mcal{N} \}$ for $0 \leq j \leq m - 2$ and $\se \in (- \iy, \se_{j+1}^{\prime}] \cap \Se$. \eeL

\bpf

By (\ref{lemA8}) of Lemma \ref{bridge} and the definition of the stopping and decision rules, \bee \{ \tx{$\mscr{H}_\ell$ with some $\ell > j$
is accepted} \} & \subseteq & \bigcup_{\ell > j } \li \{ \Up_n (\bs{\mcal{X}}_n;
\se_i^{\prime}, \se_i^{\prime \prime}) \geq \f{1}{\al_i}, \;  1 \leq i \leq \ell \; \tx{for some} \; n \in \mcal{N} \ri \}\\
& \subseteq & \li \{ \Up_n (\bs{\mcal{X}}_n;  \se_{j+1}^{\prime}, \se_{j+1}^{\prime \prime}) \geq \f{1}{\al_{j+1}} \; \tx{for some} \; n \in \mcal{N} \ri \}\\
& \subseteq & \li \{ \Up_n (\bs{\mcal{X}}_n;  \se, \se_{j+1}^{\prime \prime} ) \geq \f{1}{\al_{j+1}} \; \tx{for some} \; n \in \mcal{N} \ri \}
\eee for $0 \leq j \leq m - 2$ and $\se \in (- \iy, \se_{j+1}^{\prime}] \cap \Se$.

\epf

\beL \la{lemba}

$\{ \tx{$\mscr{H}_\ell$  with some $\ell < j$ is accepted} \}  \subseteq  \{ \Up_n (\bs{\mcal{X}}_n; \se_j^{\prime}, \se) \leq \ba_j \; \tx{for
some} \; n \in \mcal{N} \}$ for $1 \leq j \leq m - 1$ and $\se \in [\se_j^{\prime \prime}, \iy) \cap \Se$. \eeL

\bpf By (\ref{lemB8}) of Lemma \ref{bridge} and the definition of the stopping and decision rules, \bee \{ \tx{$\mscr{H}_\ell$  with some $\ell
< j$ is accepted} \} & \subseteq  & \bigcup_{\ell < j } \{ \Up_n (\bs{\mcal{X}}_n;
\se_i^{\prime}, \se_i^{\prime \prime}) \leq \ba_i, \; \ell < i < m \; \tx{for some} \; n \in \mcal{N} \}\\
& \subseteq & \{ \Up_n (\bs{\mcal{X}}_n; \se_j^{\prime}, \se_j^{\prime \prime}) \leq \ba_j \; \tx{for some} \; n \in \mcal{N} \}\\
& \subseteq & \{ \Up_n (\bs{\mcal{X}}_n; \se_j^{\prime}, \se) \leq \ba_j \; \tx{for some} \; n \in \mcal{N} \} \eee for $1 \leq j \leq m - 1$
and $\se \in [\se_j^{\prime \prime}, \iy) \cap \Se$.

\epf

\beL \la{lem88} Let $0 < j < m$ and $\se \in (\se_j^{\prime \prime}, \iy) \cap \Se$. Then, $\Pr \li \{ \Up_n (\bs{\mcal{X}}_n; \se_i^{\prime},
\se_i^{\prime \prime}) \geq \f{1}{\al_i} \; \tx{for}  \; 0 < i \leq j \mid \se \ri \} \to 1$ as the sample number $n$ tends to $N^*$. \eeL

\bpf  Let $\ba \in (0, 1)$.  By (\ref{lemB8}) of Lemma \ref{bridge},  for $0 < j < m$ and $\se \in (\se_j^{\prime \prime}, \iy) \cap \Se$, \bee
& & \Pr \li \{ \tx{There exists some $i$ such that} \; 0 < i \leq j \; \tx{and that} \; \Up_n (\bs{\mcal{X}}_n; \se_i^{\prime}, \se_i^{\prime
\prime}) < \f{1}{\al_i} \mid \se \ri \}\\
&  & \leq \sum_{i = 1}^j \Pr \li \{ \Up_n (\bs{\mcal{X}}_n; \se_i^{\prime}, \se_i^{\prime \prime}) <
\f{1}{\al_i}  \mid \se \ri \}\\
&  & \leq \sum_{i = 1}^j \li [ \Pr \li \{ \ba < \Up_n (\bs{\mcal{X}}_n; \se_i^{\prime}, \se_i^{\prime \prime}) <
\f{1}{\al_i} \mid \se \ri \}  +  \Pr \li \{ \Up_n (\bs{\mcal{X}}_n; \se_i^{\prime}, \se_i^{\prime \prime}) \leq \ba  \mid \se \ri \} \ri ]\\
&  & \leq \sum_{i = 1}^j \li [ \Pr \li \{ \ba < \Up_n (\bs{\mcal{X}}_n; \se_i^{\prime}, \se_i^{\prime \prime}) < \f{1}{\al_i} \mid
\se \ri \}  +  \Pr \li \{ \Up_n (\bs{\mcal{X}}_n; \se_i^{\prime}, \se) \leq \ba  \mid \se \ri \} \ri ]\\
&  & \leq \sum_{i = 1}^j \li [ \Pr \li \{ \ba < \Up_n (\bs{\mcal{X}}_n; \se_i^{\prime}, \se_i^{\prime \prime}) < \f{1}{\al_i} \mid \se \ri \} +
\ba  \ri ] \to j \ba
 \eee
as the sample number $n$ tends to $N^*$. But this holds for arbitrarily small $\ba \in (0, 1)$.   \epf

\beL \la{lem77} Let $0 < j < m$ and $\se \in (- \iy, \se_j^\prime) \cap \Se$.  Then, $\Pr \li \{ \Up_n (\bs{\mcal{X}}_n; \se_i^{\prime},
\se_i^{\prime \prime}) \leq \ba_i \; \tx{for} \; j \leq i < m \mid \se \ri \} \to 1$ as the sample number $n$ tends to $N^*$.

\eeL

\bpf Let $\al \in (0, 1)$.  By (\ref{lemA8}) of Lemma \ref{bridge}, for $0 < j < m$ and $\se \in (- \iy, \se_j^\prime) \cap \Se$, \bee & & \Pr
\li \{ \tx{There exists some $i$ such that $\; j \leq i < m$ and that} \; \Up_n (\bs{\mcal{X}}_n; \se_i^{\prime},
\se_i^{\prime \prime}) > \ba_i \mid \se \ri \}\\
&  & \leq \sum_{i = j }^{m-1} \Pr \li \{ \Up_n (\bs{\mcal{X}}_n; \se_i^{\prime}, \se_i^{\prime \prime}) >
\ba_i \mid \se \ri \}\\
&  & \leq \sum_{i = j }^{m-1} \li [ \Pr \li \{ \ba_i < \Up_n (\bs{\mcal{X}}_n; \se_i^{\prime}, \se_i^{\prime \prime})
< \f{1}{\al} \mid \se \ri \}  + \Pr \li \{ \Up_n (\bs{\mcal{X}}_n; \se_i^{\prime}, \se_i^{\prime \prime}) \geq \f{1}{\al} \mid \se \ri \} \ri ]\\
&  & \leq \sum_{i = j }^{m-1} \li [ \Pr \li \{ \ba_i < \Up_n (\bs{\mcal{X}}_n; \se_i^{\prime}, \se_i^{\prime \prime}) <
\f{1}{\al} \mid \se \ri \}  + \Pr \li \{ \Up_n (\bs{\mcal{X}}_n; \se, \se_i^{\prime \prime}) \geq \f{1}{\al} \mid \se \ri \} \ri ]\\
&  & \leq \sum_{i = j }^{m-1} \li [ \Pr \li \{ \ba_i < \Up_n (\bs{\mcal{X}}_n; \se_i^{\prime}, \se_i^{\prime \prime}) < \f{1}{\al} \mid \se \ri
\}  + \al \ri ]  \to (m - j) \al \eee  as the sample number $n$ tends to $N^*$. But this holds for arbitrarily small $\al \in (0, 1)$. \epf

We are now in a position to prove the theorem.

\subsection{Proof of Statements (I)--(III) }

We shall show statements (I)--(III) based on the assumption that the likelihood function $f_n (\bs{\mcal{X}}_n; \se)$ is unimodal with respect
to $\se \in \Se$ and that the sampling process will eventually terminate according to the stopping rule.

Statement (I) can be shown as follows.  Invoking  Lemmas \ref{Villeinq} and \ref{lemal}, we have \be \la{qq18}  \Pr \{ \tx{$\mscr{H}_i$ with
some $i \geq j$ is accepted} \mid \se \} \leq \Pr \li \{ \Up_n (\bs{\mcal{X}}_n;  \se, \se_j^{\prime \prime} ) \geq \f{1}{\al_j} \; \tx{for
some} \; n  \in \mcal{N} \mid \se \ri \} \leq \al_j \ee for $j = 1, \cd, m - 1$ and $\se \in (- \iy, \; \se_j^\prime ) \cap \Se$.  Making use of
Lemmas \ref{Villeinq} and \ref{lemba}, we have \be \la{qqb}
 \Pr \{ \tx{$\mscr{H}_i$  with some $i < j$ is accepted} \mid \se \} \leq \Pr \{ \Up_n (\bs{\mcal{X}}_n;  \se_j^{\prime}, \se)
\leq \ba_{m - 1} \; \tx{for some} \; n  \in \mcal{N} \mid \se \} \leq \ba_{m-1} \ee for $j = 1, \cd, m - 1$ and $\se \in   ( \se_j^{\prime
\prime}, \; \iy) \cap \Se$.  Therefore,  \bee &  & \Pr \{ \tx{Reject} \; \mscr{H}_0 \mid \se \}  =  \Pr \{ \tx{$\mscr{H}_i$ with some $i \geq 1$
is accepted} \mid \se \} \leq \al_1 \qu \tx{for $\se \in (- \iy, \; \se_1^\prime ) \cap \Se$}, \\
&  & \Pr \{ \tx{Reject} \; \mscr{H}_{m-1} \mid \se \} = \Pr \{ \tx{$\mscr{H}_i$  with some $i < {m - 1}$ is accepted} \mid \se \} \leq \ba_{m-1}
\qu  \tx{for $\se \in   ( \se_{m-1}^{\prime \prime}, \; \iy) \cap \Se$} \eee and  \bee \Pr \{ \tx{Reject} \; \mscr{H}_j \mid \se \} & = & \Pr \{
\tx{$\mscr{H}_i$ with some $i > j$ is accepted} \mid \se \} + \Pr \{ \tx{$\mscr{H}_i$  with some $i < j$ is accepted} \mid \se \} \\
& \leq & \al_{j + 1} + \ba_j \eee for $0 <  j \leq m - 2$ and $\se \in ( \se_j^{\prime \prime}, \; \se_{j + 1}^\prime ) \cap \Se$.  This proves
statement (I).

To show statement (II), let $\mbf{n}$ denote the sample number at the termination of the sampling process. Since the likelihood function $f_n
(\bs{\mcal{X}}_n; \se)$ is unimodal with respect to $\se \in \Se$, we have that for every value $n$ in the support of $\mbf{n}$, there exists an
estimator $\wh{\se}_n$, defined in terms of $\bs{\mcal{X}}_n$, such that $f_n (\bs{\mcal{X}}_n; \se)$ is nondecreasing with respect to $\se \in
\Se$ no greater than $\wh{\se}_n$ and is non-increasing with respect to $\se \in \Se$ no less than $\wh{\se}_n$. Define a sequential estimator
$\wh{\bs{\se}}$ by replacing $n$ with $\mbf{n}$, that is $\wh{\bs{\se}} = \wh{\se}_{\mbf{n}}$.  Then, $\wh{\bs{\se}}$ is a ULE of $\se$.  By the
definition of the stopping and decision rules, we have, for $j = 1, \cd, m - 1$ and every $n$ in the support of $\mbf{n}$, {\small \bee \{
\tx{Accept} \; \mscr{H}_i \; \tx{with some index $i$ no less than $j$} \} \cap \{ \mbf{n} = n \} \subseteq \li \{ \Up_n (\bs{\mcal{X}}_n;
\se_j^\prime, \se_j^{\prime \prime} ) \geq \f{1}{\al_j}, \; \mbf{n} = n \ri \} \subseteq \li \{ \wh{\bs{\se}} \geq \se_j^\prime, \; \mbf{n} = n
\ri \}. \eee} It follows that $\{ \tx{Accept} \; \mscr{H}_i \; \tx{with some index $i$ no less than $j$} \} \subseteq \li \{ \wh{\bs{\se}} \geq
\se_j^\prime \ri \}$ for $j = 1, \cd, m - 1$.

According to the second statement of Lemma 3 of \cite[version 32, Appendix A3, page 127]{Chenestimation}, we have that $\Pr \{ \tx{Accept} \;
\mscr{H}_i \; \tx{with some index $i$ no less than $j$} \mid \se \}$ is non-decreasing with respect to $\se \in \Se$ no greater than
$\se_j^\prime$ for $j = 1, \cd, m - 1$. This result together with the proven inequality (\ref{qq18}) complete the proof of Statement (II).
Similarly, to show statement (III), note that, for $j = 1, \cd, m - 1$ and every $n$ in the support of $\mbf{n}$, {\small \bee \{ \tx{Accept} \;
\mscr{H}_i \; \tx{with some index $i$ less than $j$} \} \cap \{ \mbf{n} = n \} \subseteq \li \{ \Up_n (\bs{\mcal{X}}_n; \se_j^\prime,
\se_j^{\prime \prime} ) \leq \ba_j, \; \mbf{n} = n \ri \} \subseteq \li \{ \wh{\bs{\se}} \leq \se_j^{\prime \prime}, \; \mbf{n} = n \ri \}.
\eee}

According to the first statement of Lemma 3 of \cite[version 32, Appendix A3, page 127]{Chenestimation}, we have that $\Pr \{ \tx{Accept} \;
\mscr{H}_i \; \tx{with some index $i$ less than $j$} \mid \se \}$ is non-increasing with respect to $\se \in \Se$ no less than $\se_j^{\prime
\prime}$ for $j = 1, \cd, m - 1$.  This result together with the proven inequality (\ref{qqb}) complete the proof of Statement (III).

\subsection{Proof of the Termination Property}

We shall show that the sampling process will eventually terminate according to the stopping rule under the assumption that the likelihood
function $f_n (\bs{\mcal{X}}_n; \se)$ is unimodal with respect to $\se \in \Se$ and that (\ref{TwoStop}) is satisfied. Note that for all $n$ and
$\se \in (- \iy, \se_1^\prime) \cap \Se$, \bee & & \Pr \{
\tx{The sampling process will eventually terminate according to the stopping rule} \mid \se \}\\
&  & \geq \Pr \li \{ \Up_n (\bs{\mcal{X}}_n; \se_i^{\prime}, \se_i^{\prime \prime}) \leq \ba_i \; \tx{for} \; 0 < i < m \mid \se \ri \} \eee It
follows from Lemma \ref{lem77} that $\Pr \li \{ \Up_n (\bs{\mcal{X}}_n; \se_i^{\prime}, \se_i^{\prime \prime}) \leq \ba_i \; \tx{for} \; 0 < i <
m \mid \se \ri \} \to 1$ as the sample number $n$ tends to $N^*$.  It must be true that \[ \Pr \{ \tx{The sampling process will eventually
terminate according to the stopping rule} \mid \se \} = 1
\]
for $\se \in (- \iy, \se_1^\prime) \cap \Se$.

Note that for all $n$ and $\se \in (\se_{m-1}^{\prime \prime}, \iy) \cap \Se$, \bee &  & \Pr \{ \tx{The sampling process will eventually
terminate according to the stopping rule} \mid \se \}\\
&  & \geq \Pr \li \{ \Up_n (\bs{\mcal{X}}_n; \se_i^{\prime}, \se_i^{\prime \prime}) \geq \f{1}{\al_i} \; \tx{for} \; 0 <  i < m \mid \se \ri \}
\eee It follows from Lemma \ref{lem88} that $\Pr \li \{ \Up_n (\bs{\mcal{X}}_n; \se_i^{\prime}, \se_i^{\prime \prime}) \geq \f{1}{\al_i} \;
\tx{for} \; 0 <  i < m \mid \se \ri \} \to 1$ as the sample number $n$ tends to $N^*$.  It must be true that
\[
\Pr \{ \tx{The sampling process will eventually terminate according to the stopping rule} \mid \se \} = 1 \]
 for $\se \in (\se_{m-1}^{\prime
\prime}, \iy) \cap \Se$.  By Lemmas \ref{lem88}, \ref{lem77} and Bonferroni's inequality,  we have \bee \Pr \li \{ \Up_n (\bs{\mcal{X}}_n;
\se_i^{\prime}, \se_i^{\prime \prime}) \geq \f{1}{\al_i},  \; 0 <  i \leq j \; \tx{and} \; \Up_n (\bs{\mcal{X}}_n; \se_i^{\prime}, \se_i^{\prime
\prime}) \leq \ba_i,  \; j < i < m \mid \se \ri \}  \to 1 \eee for $j = 1, \cd, m - 2$ and $\se \in (\se_j^{\prime \prime}, \se_{j + 1}^\prime)
\cap \Se$, as the sample number $n$ tends to $N^*$.  Note that for all $n$, $j = 1, \cd, m - 2$ and $\se \in (\se_j^{\prime \prime}, \se_{j +
1}^\prime) \cap \Se$, \bee
&  & \Pr \{ \tx{The sampling process will eventually terminate according to the stopping rule} \mid \se \}\\
&  & \geq \Pr \li \{ \Up_n (\bs{\mcal{X}}_n; \se_i^{\prime}, \se_i^{\prime \prime}) \geq \f{1}{\al_i},  \; 0 <  i \leq j \; \tx{and} \; \Up_n
(\bs{\mcal{X}}_n; \se_i^{\prime}, \se_i^{\prime \prime}) \leq \ba_i,  \; j < i < m \mid \se \ri \}. \eee Thus, it must be true that
\[
\Pr \{ \tx{The sampling process will eventually terminate according to the stopping rule} \mid \se \} = 1 \] for $j = 1, \cd, m - 2$ and $\se
\in (\se_j^{\prime \prime}, \se_{j + 1}^\prime) \cap \Se$.

It remains to show that the sampling process will eventually terminate according to the stopping rule for $\se \in [\se_j^{\prime},
\se_j^{\prime \prime}] \cap \Se$ with $j = 1, \cd, m - 1$.  By Lemma \ref{lem88}, for $1 < j < m$ and all $\se \in [\se_j^{\prime},
\se_j^{\prime \prime}] \cap \Se \subseteq (\se_{j-1}^{\prime \prime}, \iy) \cap \Se$, \be \la{qq1} \Pr \li \{ \Up_n (\bs{\mcal{X}}_n;
\se_i^{\prime}, \se_i^{\prime \prime}) \geq \f{1}{\al_i} \; \tx{for} \; 0 <  i < j \mid \se \ri \} \to 1 \ee as the sample number $n$ tends to
$N^*$.  By Lemma \ref{lem77}, for $0 \leq j < m - 1$ and $\se \in [\se_j^{\prime}, \se_j^{\prime \prime}] \cap \Se \subseteq (- \iy,
\se_{j+1}^{\prime} ) \cap \Se$,  \be \la{qq2}
 \Pr \li \{ \Up_n (\bs{\mcal{X}}_n; \se_i^{\prime}, \se_i^{\prime \prime}) \leq \ba_i \; \tx{for} \; j < i < m \mid
\se \ri \} \to 1 \ee as the sample number $n$ tends to $N^*$.  By the assumption associated with (\ref{TwoStop}), for $j = 1, \cd, m - 1$ and
$\se \in [\se_j^{\prime}, \se_j^{\prime \prime}] \cap \Se$, \be \la{qq3} \Pr \li \{ \Up_n (\bs{\mcal{X}}_n; \se_j^{\prime}, \se_j^{\prime
\prime}) \geq \f{1}{\al_j} \; \tx{or} \; \Up_n (\bs{\mcal{X}}_n; \se_j^{\prime}, \se_j^{\prime \prime}) \leq \ba_j \mid \se \ri \} \to 1 \ee as
the sample number $n$ tends to
$N^*$.  Note that \bee &  & \{ \tx{The sampling process will eventually terminate according to the stopping rule}  \}\\
 &   & \supseteq \li \{ \Up_n
(\bs{\mcal{X}}_n; \se_i^{\prime}, \se_i^{\prime \prime}) \geq \f{1}{\al_i} \; \tx{for} \; 0 <  i \leq j  \; \tx{and} \;
\Up_n (\bs{\mcal{X}}_n; \se_i^{\prime}, \se_i^{\prime \prime}) \leq \ba_i \; \tx{for} \; j < i < m \ri \}\\
&  & \qu \bigcup \li \{ \Up_n (\bs{\mcal{X}}_n; \se_i^{\prime}, \se_i^{\prime \prime}) \geq \f{1}{\al_i} \; \tx{for} \; 0 <  i < j  \; \tx{and}
\; \Up_n (\bs{\mcal{X}}_n; \se_i^{\prime}, \se_i^{\prime \prime}) \leq \ba_i \; \tx{for} \; j \leq i < m \ri \}\\
 &   & = \li \{ \Up_n (\bs{\mcal{X}}_n;
\se_i^{\prime}, \se_i^{\prime \prime}) \geq \f{1}{\al_i} \; \tx{for} \; 0 <  i < j  \ri \} \bigcap \li \{ \Up_n (\bs{\mcal{X}}_n;
\se_i^{\prime}, \se_i^{\prime \prime})
\leq \ba_i \; \tx{for} \; j < i < m \ri \}\\
&  & \qu \bigcap \li \{ \Up_n (\bs{\mcal{X}}_n; \se_j^{\prime}, \se_j^{\prime \prime}) \geq \f{1}{\al_j} \; \tx{or} \; \Up_n (\bs{\mcal{X}}_n;
\se_j^{\prime}, \se_j^{\prime \prime}) \leq \ba_j \ri \} \eee for $j = 2, \cd, m - 2$. Making use of this observation, (\ref{qq1}), (\ref{qq2}),
(\ref{qq3}) and Bonferroni's inequality,  we have {\small \bee & & \Pr \{ \tx{The sampling process will
eventually terminate according to the stopping rule} \mid \se \}\\
&  & \geq \Pr \li \{ \Up_n (\bs{\mcal{X}}_n; \se_i^{\prime}, \se_i^{\prime \prime}) \geq \f{1}{\al_i} \; \tx{for} \; 0 <  i < j  \mid \se \ri \}
+
\Pr \li \{ \Up_n (\bs{\mcal{X}}_n; \se_i^{\prime}, \se_i^{\prime \prime}) \leq \ba_i \; \tx{for} \; j < i < m \mid \se \ri \}\\
&  & \; \; \; \; + \Pr \li \{ \Up_n (\bs{\mcal{X}}_n; \se_j^{\prime}, \se_j^{\prime \prime}) \geq \f{1}{\al_j} \; \tx{or} \;
\Up_n (\bs{\mcal{X}}_n; \se_j^{\prime}, \se_j^{\prime \prime}) \leq \ba_j  \mid \se \ri \} - 3\\
&  &  \to 1 \eee} for $\se \in [\se_j^{\prime}, \se_j^{\prime \prime}] \cap \Se$ with $j = 2, \cd, m - 2$, as the sample number $n$ tends to
$N^*$.   By Bonferroni's inequality,  we have {\small \bee &  & \Pr \{ \tx{The sampling process will eventually terminate according to the stopping rule} \mid \se \}\\
 &   & \geq \Pr \li \{ \Up_n (\bs{\mcal{X}}_n; \se_1^{\prime}, \se_1^{\prime \prime}) \geq \f{1}{\al_1} \; \tx{or} \;  \Up_n (\bs{\mcal{X}}_n;
\se_1^{\prime}, \se_1^{\prime \prime}) \leq \ba_1  \mid \se \ri \} + \Pr \li \{ \Up_n
(\bs{\mcal{X}}_n; \se_i^{\prime}, \se_i^{\prime \prime}) \leq \ba_i,  \; 1 < i < m \mid \se \ri \} - 2\\
&  & \to  1 \eee} for $\se \in [\se_1^{\prime}, \se_1^{\prime \prime}] \cap \Se$, as the sample number $n$ tends to $N^*$.  Again by
Bonferroni's inequality,  we have {\small \bee & & \Pr \{ \tx{The sampling process will
eventually terminate according to the stopping rule} \mid \se \}\\
&  & \geq \Pr \li \{ \Up_n (\bs{\mcal{X}}_n; \se_{m-1}^{\prime}, \se_{m-1}^{\prime \prime}) \geq \f{1}{\al_{m-1}} \; \tx{or} \; \Up_n
(\bs{\mcal{X}}_n; \se_{m-1}^{\prime}, \se_{m-1}^{\prime \prime}) \leq \ba_{m-1}  \mid \se \ri \} \\
&  & \; \; \; \; + \Pr \li \{ \Up_n (\bs{\mcal{X}}_n; \se_i^{\prime}, \se_i^{\prime \prime})
\geq \f{1}{\al_i},  \; 0 <  i < m-1 \mid \se \ri \} - 2 \\
&  &  \to 1 \eee} for $\se \in [\se_{m-1}^{\prime}, \se_{m-1}^{\prime \prime}] \cap \Se$, as the sample number $n$ tends to $N^*$.  Therefore,
we have shown that, with probability $1$,  the sampling process will eventually terminate according to the stopping rule for $\se \in
[\se_j^{\prime}, \se_j^{\prime \prime}]$ for $j = 1, \cd, m - 1$.  This completes the proof of the theorem.

\section{Proof of Theorem \ref{MSPRT_Simple}} \la{MSPRT_Simple_app}

We need some preliminary results.

\beL \la{cooll}

For $0 \leq j < m - 1$,  \be \la{coola} \{ \tx{$\mscr{H}_\ell$ with some $\ell
> j$ is accepted} \} \subseteq \li \{ \Up_n (\bs{\mcal{X}}_n; \se_j, \se_{j + 1}) \geq \f{1}{\al_{j + 1}}\; \tx{for some} \; n \in \mcal{N} \ri \}
\ee Similarly, \be \la{coolb} \{ \tx{$\mscr{H}_\ell$  with some $\ell < j$ is accepted} \}  \subseteq  \li \{ \Up_n (\bs{\mcal{X}}_n; \se_{j -
1}, \se_j) \leq \ba_j \; \tx{for some} \; n  \in \mcal{N} \ri \} \ee for $1 \leq j < m$. \eeL

\bpf By the definition of the stopping and decision rules, \bee \{ \tx{$\mscr{H}_\ell$ with some $\ell
> j$ is accepted} \} & \subseteq & \bigcup_{\ell > j } \li \{ \Up_n (\bs{\mcal{X}}_n; \se_{i - 1}, \se_i)
\geq \f{1}{\al_i}, \;  1 \leq i \leq \ell \; \tx{for some} \; n \in \mcal{N} \ri \}\\
& \subseteq & \li \{ \Up_n (\bs{\mcal{X}}_n; \se_j, \se_{j + 1}) \geq \f{1}{\al_{j + 1}}\; \tx{for some} \; n \in \mcal{N} \ri \} \eee for $0
\leq j < m - 1$.  Similarly, by the definition of the stopping and decision rules, \bee \{ \tx{$\mscr{H}_\ell$  with some $\ell < j$ is
accepted} \} & \subseteq & \bigcup_{\ell < j } \li \{ \Up_n (\bs{\mcal{X}}_n; \se_{i - 1}, \se_i) \leq \ba_i, \; \ell < i \leq m - 1 \; \tx{for
some} \; n \in \mcal{N} \ri \}\\
& \subseteq & \li \{ \Up_n (\bs{\mcal{X}}_n; \se_{j - 1}, \se_j) \leq \ba_j \; \tx{for some} \; n \in \mcal{N} \ri \} \eee for $1 \leq j < m$.

\epf

We are now in a position to prove the theorem.  It follows from (\ref{coola}) of Lemma \ref{cooll} that \bee \Pr \{ \tx{Reject} \; \mscr{H}_0
\mid \se_0 \} & = & \Pr \{ \tx{$\mscr{H}_\ell$ with some $\ell
> 0$ is accepted} \mid \se_0 \}\\
& \leq & \Pr \li \{ \Up_n (\bs{\mcal{X}}_n; \se_0, \se_1) \geq \f{1}{\al_1} \; \tx{for some} \; n \in \mcal{N} \mid \se_0 \ri \} \leq \al_1.
\eee It follows from (\ref{coolb}) of Lemma \ref{cooll} that \bee \Pr \{ \tx{Reject} \; \mscr{H}_{m-1} \mid \se_{m-1} \} & = & \Pr
\{ \tx{$\mscr{H}_\ell$  with some $\ell < {m - 1}$ is accepted} \mid \se_{m-1} \}\\
& \leq & \Pr \li \{ \Up_n (\bs{\mcal{X}}_n; \se_{m - 2}, \se_{m-1}) \leq \ba_{m-1} \; \tx{for some} \; n  \in \mcal{N} \mid \se_{m-1} \ri \}
\leq \ba_{m-1}.
\eee It follows from (\ref{coola}) and (\ref{coolb}) of Lemma \ref{cooll} that {\small \bee &  & \Pr \{ \tx{Reject} \; \mscr{H}_j \mid \se_j \}\\
 &   & = \Pr \{ \tx{$\mscr{H}_\ell$  with some
$\ell > j$ is accepted} \mid \se_j \} + \Pr \{ \tx{$\mscr{H}_\ell$  with some $\ell < j$ is accepted} \mid \se_j \}  \\
&  & \leq  \Pr \li \{ \Up_n (\bs{\mcal{X}}_n; \se_j, \se_{j + 1}) \geq \f{1}{\al_{j + 1}}
 \; \tx{for some} \; n \in \mcal{N} \mid \se_j \ri \} + \Pr \li \{ \Up_n (\bs{\mcal{X}}_n; \se_{j-1}, \se_j) \leq \ba_j
  \; \tx{for some} \; n  \in \mcal{N} \mid \se_j \ri \}\\
&  & \leq   \al_{j+1} + \ba_{j} \eee} for $1 \leq j \leq m - 2$.

To show that the sampling process will eventually terminate according to the stopping rule with probability $1$, it suffices to apply the
argument of the proof of the termination property of Theorem \ref{MSPRT_Composite} in Appendix \ref{MSPRT_Composite_app} to the following
hypotheses
\[
\mcal{H}_0: \se \leq \vse_1, \qu \mcal{H}_1: \vse_1 < \se \leq \vse_2, \qu \ldots, \qu \mcal{H}_{m-2}: \vse_{m-2} < \se \leq \vse_{m-1}, \qu
\mcal{H}_{m-1}: \se > \vse_{m-1}
\]
with $\vse_i = \f{ \se_{i - 1} + \se_i} {2}, \; i = 1, \cd, m-1$ and indifference zone $\cup_{i=1}^{m-1} (\se_{i-1}, \se_i)$.   This concludes
the proof of the theorem.

\section{Proof of Theorem \ref{genST}} \la{genST_app}

For simplicity of notations, define $Y = \ln \f{ f(X; \se^{\prime \prime})   }{ f(X; \se^{\prime}) }$. Let $\mu$ and $\nu$ denote, respectively,
the mean and variance of $Y$ associated with $\se \in \Se$. Let $Y_1, Y_2, \cd$ be i.i.d. samples of $Y$. Define $Z_n = \f{ \sum_{i = 1}^n (Y_i
- \mu)  } { \sq{n \nu}  }$ for $n = 1, 2, \cd$.  By the central limit theorem, $Z_n$ converges in distribution to a Gaussian random variable,
$Z$,  with zero mean and unit variance.  Note that \bee \Pr \li \{ \ln \ba < \sum_{i = 1}^n Y_i < \ln \f{1}{\al} \mid \se \ri \}  =  \Pr \li \{
\f{ \ln \ba - n \mu } { \sq{n \nu} } < Z_n < \f{ \ln \f{1}{\al} - n \mu } { \sq{n \nu}  } \mid \se \ri \}, \qqu n = 1, 2, \cd \eee for $\se \in
\Se$.

In the case of $\mu > 0$, we have \bee \Pr \li \{  \ln \ba < \sum_{i = 1}^n Y_i < \ln \f{1}{\al} \mid \se \ri \}  \leq  \Pr \li \{  Z_n < \f{
\ln \f{1}{\al} - n \mu  } { \sq{n \nu}  } \mid \se \ri \} \to 0  \eee for $\se \in \Se$ as $n \to \iy$.  To show this, let $\vep > 0$.  Let $z$
be a number such that $\Pr \{ Z < z \} < \f{\vep}{2}$. Let $n$ be chosen such that $z > \f{ \ln \f{1}{\al} - n \mu } { \sq{n \nu} }$ and that $|
\Pr \{ Z_n < z \} - \Pr \{ Z < z \} | < \f{\vep}{2}$.  By the triangle inequality, \bee \Pr \li \{ Z_n < \f{ \ln \f{1}{\al} - n \mu  } { \sq{n
\nu}  } \mid \se \ri \}  \leq \Pr \li \{  Z_n < z  \mid \se \ri \} \leq \Pr \{ Z < z \}  + | \Pr \{ Z_n < z \} - \Pr \{ Z < z \} | < \vep. \eee

In the case of $\mu = 0$, we have \bee \Pr \li \{  \ln \ba < \sum_{i = 1}^n Y_i < \ln \f{1}{\al} \mid \se \ri \}  =  \Pr \li \{ \f{ \ln \ba } {
\sq{n \nu} } < Z_n < \f{ \ln \f{1}{\al} } { \sq{n \nu}  } \mid \se \ri \} \to 0 \eee for $\se \in \Se$ as $n \to \iy$.

In the case of $\mu < 0$, we have \bee \Pr \li \{  \ln \ba < \sum_{i = 1}^n Z_i < \ln \f{1}{\al} \mid \se \ri \}  \leq  \Pr \li \{  Y_n > \f{
\ln \al - n \mu  } { \sq{n \nu}  } \mid \se \ri \} \to 0  \eee for $\se \in \Se$ as $n \to \iy$.

This completes the proof of the theorem.

\section{Proof of Theorem \ref{STChe}} \la{STChe_app}

As in the proof of Theorem \ref{genST} in Appendix \ref{genST_app}, for simplicity of notations, define $Y = \ln \f{ f(X; \se^{\prime \prime})
}{ f(X; \se^{\prime}) }$. Let $\mu$ and $\nu$ denote, respectively, the mean and variance of $Y$ associated with $\se \in \Se$.  By the
assumption of the theorem, we have $\mu > 0$ and $0 < \nu < \iy$.  Let $Y_1, Y_2, \cd$ be i.i.d. sample of $Y$.   Note that \bee \Pr \li \{
\Up_n (\bs{\mcal{X}}_n; \se^{\prime}, \se^{\prime \prime}) < \f{1}{\al} \mid \se \ri
\} & = & \Pr \li \{ \sum_{i=1}^n Y_i  \geq \ln \f{1}{\al} \mid \se \ri \}\\
& = & \Pr \li \{ \f{ \sum_{i=1}^n ( Y_i - \mu ) }{n} < \f{ \ln \f{1}{\al} }{n} - \mu \mid \se \ri
\}\\
& \leq & \Pr \li \{ \f{ \sum_{i=1}^n ( Y_i - \mu ) }{n} < - \f{\mu}{2} \mid \se \ri \}
 \eee
for $n > \f{2 \ln \f{1}{\al} }{\mu}$. By Chebyshev's inequality,
\[
\Pr \li \{ \Up_n (\bs{\mcal{X}}_n; \se^{\prime}, \se^{\prime \prime}) < \f{1}{\al} \mid \se \ri \}  \leq  \Pr \li \{ \li | \f{ \sum_{i=1}^n (Y_i
- \mu ) }{n} \ri | > \f{\mu}{2} \mid \se \ri \} \leq \f{ 4 \nu } { n \mu^2 } \to 0
\]
as $n \to \iy$.  This establishes statement (I).  In a similar manner, we can show statement (II).  This completes the proof of the theorem.

\sect{Proof of Theorem \ref{expST} }  \la{expST_app}

For simplicity of notations, define $Y = T(X)$.  We need some preliminary results.

\beL \la{derlema} The derivative of $\exp \li ( u(\se) z - v (\se) \ri )$ with respect to $\se$ is equal to $(z - \se) \exp \li ( u(\se) z - v
(\se) \ri )  \f{d u(\se)}{d \se}$.  \eeL

\bpf Since $\f{d v (\se)}{d \se} = \se \f{d u(\se)}{d \se}$ for $\se \in \Se$,  by the chain rule of differentiation, we have that the
derivative of $\exp \li ( u(\se) z - v (\se) \ri )$ with respect to $\se$ is equal to $(z - \se) \exp \li ( u(\se) z - v (\se) \ri ) \f{d
u(\se)}{d \se}$.

\epf

 \beL  \la{meanl} The expectation of $Y$ is equal to $\se$.
 \eeL

 \bpf Let $\psi (.)$ be the inverse function of $u(.)$ such that $u(\psi (\ze) ) = \ze$ for $\ze \in \{ u(\se): \se \in \Se \}$. Define
compound function $w(.)$ such that $w (\zeta) = v (\psi (\zeta))$ for $\ze \in \{ u(\se): \se \in \Se \}$. For simplicity of notations, we
abbreviate $\psi (\zeta)$ as $\psi$ when this can be done without causing confusion.  Putting $\zeta = u(\se)$, we have {\small \bee  & & \bb{E}
\li [ \exp ( t Y ) \ri ] = \bb{E} \li [ \exp \li ( t T ( X ) \ri ) \ri ] =   \int h (x)  \exp \li ( ( \zeta +
t ) T ( x ) -  w (\zeta) \ri ) dx\\
& = &  \exp \li (  w (\zeta + t) - w (\zeta) \ri ) \int  h (x) \exp \li ( ( \zeta + t ) T ( x ) - w (\zeta + t) \ri ) dx  =   \exp \li (  w
(\zeta + t) - w (\zeta) \ri ). \eee} By the defining relationship $u(\psi (\ze) ) = \ze$, the assumption that $\f{ d v ( \se ) }{ d \se } = \se
\f{ d u( \se ) }{d \se}$,
 and the chain rule of differentiation, we have \be
\la{use8} \f{ d w (\zeta ) }{d \zeta} = \f{ d v (\psi) }{d \psi} \f{ d \psi }{ d \zeta} = \psi \f{ d u(\psi) }{d \psi} \f{ d \psi }{ d \zeta} =
\psi \f{ d u(\psi) }{d \ze} = \psi \f{ d \ze }{d \ze} = \psi (\ze). \ee  By virtue of (\ref{use8}), the derivative of $w (\zeta + t) - w
(\zeta)$ with respect to $t$ is given by $\f{d w (\zeta + t) }{d t} = \psi ( \zeta + t)$,  which is equal to $\psi ( \zeta ) = \se$ for $t = 0$.
Thus, $\bb{E} [ Y ] = \se$, which implies that the sample mean of $Y$ is also an unbiased estimator of $\se$.

\epf

\beL  \la{varl} The variance of $Y$ is equal to $\f{1}{ \f{d u(\se)}{d \se} }$.

\eeL

\bpf Now we are in a position to compute the variance of $Y$.  Recall that
\[
\f{d \bb{E} \li [ \exp ( t Y  ) \ri ] }{d t} = \f{d w (\zeta + t) }{d (\zeta + t)} \exp \li ( w (\zeta + t) - w (\zeta) \ri ) = \psi(\zeta + t)
\exp \li ( w (\zeta + t) - w (\zeta) \ri ).
\]
Hence, \bee &  & \f{d^2 \bb{E} \li [ \exp ( t Y  ) \ri ] }{d t^2} =  \psi^2(\zeta + t) \exp \li ( w (\zeta + t) - w (\zeta) \ri ) + \f{d
\psi(\zeta + t) }{d t} \exp \li ( w (\zeta + t) - w (\zeta) \ri )\\
&  & =  \psi^2(\zeta + t) \exp \li ( w (\zeta + t) - w (\zeta) \ri ) + \f{d \psi(\zeta + t) }{d (\zeta + t)} \exp \li ( w (\zeta + t) - w
(\zeta) \ri ). \eee Therefore, $\bb{E} [ Y^2 ] = \psi^2(\zeta) + \f{d \psi(\zeta) }{d \zeta}$.  To compute $\f{d  \psi(\zeta) }{d \ze}$, we
differentiate both sides of the defining relationship with respect to $\ze$ to obtain $\f{d u}{d \psi} \f{d \psi} {d \ze} = 1$, which implies
that $\f{d \psi} {d \ze} = \f{1}{ \f{d u}{d \psi} } = \f{1}{ \f{d u(\se)}{d \se} }$,  where we have used $\se = \psi (\ze)$ to obtain the last
equality.  Therefore, $\bb{E} [ Y^2 ] = \psi^2(\zeta) + \f{1}{ \f{d u(\se)}{d \se} } = \se^2 + \f{1}{ \f{d u(\se)}{d \se} }$,  which implies
that
\[
\mrm{Var} [Y ] = \bb{E} [ Y^2 ] - \bb{E}^2 [ Y] = \se^2 + \f{1}{ \f{d u(\se)}{d \se} } - \se^2 = \f{1}{ \f{d u(\se)}{d \se} }.
\]
\epf

We are now in a position to prove the theorem.  Since $\f{d u(\se)}{d \se}
> 0$ for $\se \in \Se$, from Lemma \ref{derlema}, we have that the derivative of $\exp \li ( u(\se) z - v (\se) \ri
)$ with respect to $\se$ is positive for $\se < z$ and negative for $\se > z$.  This implies that $\exp \li ( u(\se) z - v (\se) \ri )$ is
monotonically increasing with respect to $\se$ less than $z$ and monotonically decreasing with respect to $\se$ greater than $z$.  Since {\small
$f_n(\bs{\mcal{X}}_n; \se) =  \li [ \exp \li ( u(\se) \f{ \sum_{i=1}^n T (X_i) }{n} - v (\se) \ri ) \ri ]^n \prod_{i = 1}^n h(X_i)$}, it follows
that  $f_n(\bs{\mcal{X}}_n; \se)$ is unimodal with respect to $\se \in \Se$.

Let $X_1, X_2, \cd$ be i.i.d. samples of $X$.  For parameter values $\se^{\prime}, \; \se^{\prime \prime} \in \Se$ with $\se^{\prime} <
\se^{\prime \prime}$, the likelihood ratio is
\[
\Up_n ( \bs{\mcal{X}}_n; \se^{\prime}, \se^{\prime \prime} )  = \f{ \exp [ u (\se^{\prime \prime}) \sum_{i = 1}^n T(X_i) - n v (\se^{\prime
\prime} ) ] } { \exp [ u (\se^{\prime}) \sum_{i = 1}^n T(X_i) - n v (\se^{\prime} )] }.
\]
Note that for $n = 1, 2, \cd$,  \bee \Pr \li \{  \ba < \Up_n ( \bs{\mcal{X}}_n; \se^{\prime}, \se^{\prime \prime} ) < \f{1}{\al} \ri \} & = &
\Pr \li \{ \f{ n [ v(\se^{\prime \prime} ) -  v(\se^{\prime} ) ] +  \ln \ba }{ u(\se^{\prime \prime}) - u(\se^{\prime}) }  < \sum_{i = 1}^n
T(X_i) < \f{ n [
v(\se^{\prime \prime} ) -  v(\se^{\prime} ) ] + \ln \f{1}{\al} }{ u(\se^{\prime \prime}) - u(\se^{\prime}) } \ri \}\\
& = & \Pr \{ n \ro - a  < \sum_{i = 1}^n T(X_i)  < n \ro + b \}\\
& = &  \Pr \li \{ \f{ n (\ro - \se) - a }{ \sq{n} \si }  < Z_n  < \f{ n (\ro - \se) + b }{ \sq{n} \si } \ri \}  \eee where
\[
\ro = \f{ v(\se^{\prime \prime} ) -  v(\se^{\prime} ) } { u(\se^{\prime \prime}) - u(\se^{\prime}) }, \qu a = - \f{  \ln \ba }{ u(\se^{\prime
\prime}) - u(\se^{\prime}) }, \qu b = \f{ \ln \f{1}{\al} }{ u(\se^{\prime \prime}) - u(\se^{\prime}) }
\]
and  \[ Z_n = \f{ \sum_{i = 1}^n T(X_i) - n \se }{ \sq{n} \si }, \qqu n = 1, 2, \cd,
\]
with $\si^2 = \f{1}{ \f{d \eta (\se)}{d \se} }$ being the variance of $T(X)$. From Lemmas \ref{meanl} and \ref{varl}, we know that $T(X)$ is a
random variable with mean $\se$ and variance $\si^2$.  By the central limit theorem, $Z_n$ converges to a Gaussian random variable with zero
mean and unit variance as $n$ tends to infinity.  Consequently,
\[
\Pr \li \{ \f{ n (\ro - \se) - a }{ \sq{n} \si }  < Z_n  < \f{ n (\ro - \se) + b }{ \sq{n} \si } \ri \}  \to 0
\]
as $n \to \iy$, which can be readily shown by considering the cases of $\se > \ro, \; \se = \ro$ and $\se < \ro$ as in the proof of Theorem
\ref{genST} in Appendix \ref{genST_app}.  This completes the proof of the theorem.

\sect{Proof of Theorem \ref{CI sequence Continuous} }  \la{CI sequence Continuous_app}

We need a preliminary result.

\beL

\la{NewMart}

Suppose that $(X_n)_{n \in \bb{N}}$ is a discrete-time process parameterized by $\se \in \Se$ such that for any $n$, the conditional probability
density or mass function of $X_1, \cd, X_{n-1}$ given the value of $X_n$ does not depend on $\se$.  Let $\{ \mscr{F}_n \}$ be a natural
filtration such that for $n \in \bb{N}$, where $\mscr{F}_n$ is $\si$-algebra generated by $X_1, \cd, X_n$.  Then, for any parameter values
$\se_0$ and $\se_1$, $\li \{ \f{ f_n (X_n, \se_1) }{f_n(X_n, \se_0)} \ri \}_{n \in \bb{N}}$ is a martingale process with respect to the
filtration $\{ \mscr{F}_n \}$ and the probability measure associated with $\se_0$.

\eeL

\bpf For simplicity of notations, let $\mbf{x}_n = (x_1, \cd, x_n)$ for $n = 1, 2, \cd$.   First, consider the case that the PDF exists.  By the
assumption of the lemma, we have $\f{ f_{\bs{\mcal{X}}_n} (  \mbf{x}_n; \se_1 ) } {  f_{X_n} ( x_n; \se_1 )} = \f{ f_{\bs{\mcal{X}}_n} (
\mbf{x}_n; \se_0 ) } {  f_{X_n} ( x_n; \se_0 )}$ or equivalently, \be \la{eqkey}
 \f{ f_{\bs{\mcal{X}}_n} (  \mbf{x}_n; \se_1 ) } {  f_{\bs{\mcal{X}}_n} (  \mbf{x}_n; \se_0 ) } = \f{ f_{X_n} ( x_n; \se_1 ) } {  f_{X_n} ( x_n; \se_0 )}.
  \ee Let $d \mbf{x}_n = d x_1 \cd d x_n$ for $n = 1, 2, \cd$.
  Let $\bb{P}_{\se_0}$ denotes the probability measure associated with $\se_0 \in
  \Se$.  It follows from (\ref{eqkey}) that for arbitrary $S \subseteq \bb{R}^n$, \bee \int_{\bs{\mcal{X}}_n \in S}
  \f{ f_{X_{n+1}} ( X_{n+1}; \se_1 ) } {  f_{X_{n+1}} ( X_{n+1}; \se_0 )} d \bb{P}_{\se_0} & = &
  \int_{\mbf{x}_n \in S \atop{ x_{n+1} \in \bb{R}} } \f{ f_{X_{n+1}} ( x_{n+1}; \se_1 ) } {  f_{X_{n+1}} ( x_{n+1}; \se_0 )}
  f_{\bs{\mcal{X}}_{n+1}} (  \mbf{x}_{n+1}; \se_0 ) \; d \mbf{x}_{n+1}\\
&   = & \int_{\mbf{x}_n \in S \atop{ x_{n+1} \in \bb{R}} } \f{ f_{\bs{\mcal{X}}_{n+1}} (  \mbf{x}_{n+1}; \se_1 ) }{ f_{\bs{\mcal{X}}_{n+1}} (
\mbf{x}_{n+1}; \se_0 ) } f_{\bs{\mcal{X}}_{n+1}} (  \mbf{x}_{n+1}; \se_0 )
 \; d \mbf{x}_{n+1}\\
&   = & \int_{\mbf{x}_n \in S \atop{ x_{n+1} \in \bb{R}} } f_{\bs{\mcal{X}}_{n+1}} (  \mbf{x}_{n+1}; \se_1 ) \; d \mbf{x}_{n+1}\\
&  = & \int_{\mbf{x}_n \in S} \li [ \int_{x_{n+1} \in \bb{R}} f_{\bs{\mcal{X}}_{n+1}} (  \mbf{x}_{n+1}; \se_1 ) \; d x_{n+1} \ri ] d \mbf{x}_n\\
& = & \int_{\mbf{x}_n \in S  } f_{\bs{\mcal{X}}_{n}} (  \mbf{x}_{n}; \se_1 ) \; d \mbf{x}_n =  \int_{\mbf{x}_n \in S  } \f{
f_{\bs{\mcal{X}}_{n}} (  \mbf{x}_{n}; \se_1 )  }{  f_{\bs{\mcal{X}}_{n}} (  \mbf{x}_{n}; \se_0 ) } f_{\bs{\mcal{X}}_{n}} (
\mbf{x}_{n}; \se_0 ) \; d \mbf{x}_n\\
& = & \int_{\mbf{x}_n \in S  } \f{ f_{X_n} ( x_n; \se_1 ) } {  f_{X_n} ( x_n; \se_0 )} f_{\bs{\mcal{X}}_{n}} ( \mbf{x}_{n}; \se_0 ) \; d
\mbf{x}_n = \int_{\bs{\mcal{X}}_n \in S } \f{ f_{X_n} ( X_n; \se_1 ) } {  f_{X_n} ( X_n; \se_0 )}  d \bb{P}_{\se_0}, \eee which implies that
$\li \{ \f{ f_n (X_n, \se_1) }{f_n(X_n, \se_0)}, \; \mscr{F}_n \ri \}_{n \in \bb{N}}$ is a martingale with respect to the filtration $\{
\mscr{F}_n \}$ and the probability measure associated with $\se_0$.  In the case that the PMF exists, the integration in the above is replaced
by summation.

\epf

We are now in a position to prove the theorem.  Define $Y_t = \Up_t (X_t; \se_0, \se_1)$ for $t \in [0, \iy)$.  Define
\[
Q_k = \{0\} \cup \li \{ \f{q}{p}: \gcd(p,q) = 1; \; p, q \in \bb{N}; \; p \leq k \ri \}
\]
for $k = 1, 2, \cd$, where $\gcd(p,q)$ denotes the greatest common divider of $p$ and $q$.  Let $\wt{Q}$ denote the set of non-negative rational
numbers. Define $E_j = \{ \om \in \Om: \sup_{ t \in Q_j  } Y_t (\om)
> \f{1}{\de}  \}$.  Then,
\[
Q_j \subseteq Q_{j + 1} \Rightarrow \sup_{t \in  Q_j} Y_t (\om) \leq \sup_{t \in  Q_{j+1}} Y_t (\om) \Rightarrow E_j \subseteq E_{j + 1}.
\]
Define $E_\iy = \li \{  \om \in \Om: \sup_{ t \in \wt{Q} }  Y_t (\om) > \f{1}{\de} \ri \}$.  It is easy to show that $E_\iy = \cup_{j = 0}^\iy
E_j$. As a consequence of the continuity of the probability measure,  $\Pr \{ E_\iy \} = \lim_{n \to \iy} \Pr \{ E_n \}$.  By Lemma
\ref{NewMart}, $\{ Y_t, \; t \in Q_j \}$ is a martingale process. It follows from Doob's super-martingale inequality that $\Pr \{ E_j \} \leq
\de \; \bb{E} [ Y_0 ]$. This implies that $\Pr \{ E_\iy \} = \lim_{j \to \iy} \Pr \{ E_j \} \leq \de \; \bb{E} [ Y_0 ]$.   We claim that $\sup_{
t \in \wt{Q} } Y_t (\om) = \sup_{t \in [0, \iy)} Y_t (\om)$.  To show this claim, note that for any $t \in [0, \iy)$, there exists a sequence
$\{q_j\}_{j = 1}^\iy$ no less than $t$ such that $Y_t (\om) = \lim_{j \to \iy} Y_{q_j} (\om)$.  That is, the sample path of $Y_t$ is
right-continuous. Observing that $Y_{q_j} (\om) \leq \sup_{ t \in \wt{Q} }  Y_t (\om)$, we have $Y_t (\om) \leq \sup_{ t \in \wt{Q} } Y_t
(\om)$, which implies that $\sup_{t \in [0, \iy)} Y_t (\om) = \sup_{ t \in \wt{Q} }  Y_t (\om)$ and thus the claim is established.  This proves
(\ref{maxcon}), that is, $\Pr \{ Y_t > \f{1}{\de} \; \tx{for some} \; t \in [0, \iy) \mid \se_0 \} \leq \de$.

By the definition of the lower confidence limit, we have $\{ L_t (X_t) \leq \se_0 \} \supseteq \li \{ \Up_t (X_t; \se_1, \se_0) \geq \f{\de}{2}
\ri \}$.  This implies that  $\{ L_t (X_t) > \se_0 \} \subseteq \li \{ \Up_t (X_t; \se_1, \se_0) < \f{\de}{2} \ri \}$ and consequently,  $\Pr \{
L_t (X_t) > \se \; \tx{for some} \; t  \mid \se \}  \leq \Pr \li \{  \Up_t (X_t; \se_1, \se) < \f{\de}{2}   \; \tx{for some} \; t \mid \se \ri
\}$ for $\se \in \Se$.  It follows from the proven inequality (\ref{maxcon}) that $\Pr \{ L_t (X_t)
> \se \; \tx{for some} \; t  \mid \se \} \leq \f{\de}{2}$.

Similarly, from the definition of the upper confidence limit, we have $\{ \Up_t (X_t; \se_0, \se_1) \geq \f{\de}{2} \} \subseteq  \{ U_t (X_t)
\geq \se_1 \}$. This implies that  $\{ U_t (X_t) < \se_1 \} \subseteq \li \{ \Up_t (X_t; \se_0, \se_1) < \f{\de}{2} \ri \}$ and consequently,
$\Pr \{ U_t (X_t) < \se \; \tx{for some} \; t  \mid \se \} \leq  \Pr \li \{  \Up_t (X_t; \se_0, \se) < \f{\de}{2} \; \tx{for some} \; t \mid \se
\ri \}$ for $\se \in \Se$.   It follows from (\ref{maxcon}) that $\Pr \{ U_t (X_t) < \se \; \tx{for some} \; t \mid \se \} \leq \f{\de}{2}$. So,
by virtue of Bonferroni's inequality, we have $\Pr \{ L_t (X_t) \leq \se \leq U_t (X_t) \; \tx{for all} \; t  \mid \se \} \geq 1 - \de$. This
completes the proof of the theorem.

\sect{Proof of Theorem \ref{PoST} }  \la{PoST_app}

We need some preliminary results.

\beL

\la{PoSTL}

For arbitrary $\al, \ba \in (0, 1)$ and $\lm, \lm^{\prime}, \lm^{\prime \prime} \in (0, \iy)$ with $\lm^{\prime} < \lm^{\prime \prime}$,  \bee &
& \lim_{t \to \iy} \Pr \li \{ \ba \leq \Up_n ( \bs{\mcal{X}}_n; \lm^{\prime}, \lm^{\prime \prime} ) \leq \f{1}{\al} \mid \lm \ri \} = 0.  \eee

\eeL

 \bpf Note that \[ \Pr \li \{ \ba \leq \Up_n ( \bs{\mcal{X}}_n; \lm^{\prime}, \lm^{\prime \prime} ) \leq \f{1}{\al} \mid \lm \ri \} = \Pr \li \{ \f{
(\lm^{\prime \prime} - \lm^{\prime}) t + \ln \ba } {\ln \f{\lm^{\prime \prime}}{\lm^{\prime}} } \leq X_t \leq \f{ (\lm^{\prime \prime} -
\lm^{\prime}) t + \ln \f{1}{\al}  } {\ln \f{\lm^{\prime \prime}}{\lm^{\prime}} } \mid \lm \ri \}.
\]
Therefore, $\Pr \li \{ \ba \leq \Up_n ( \bs{\mcal{X}}_n; \lm^{\prime}, \lm^{\prime \prime} ) \leq \f{1}{\al} \mid \lm \ri \}$ can be written as
$\Pr \{ \ro t - a \leq X_t \leq \ro t + b \mid \lm \}$, where $\ro = \f{ (\lm^{\prime \prime} - \lm^{\prime})  } {\ln \f{\lm^{\prime
\prime}}{\lm^{\prime}} }$ and $a, \; b$ are some positive numbers.  Define $Y_t = \f{ X_t - \lm t }{\sq{\lm t}}$.  Then, \[ \Pr \{ \ro t - a
\leq X_t \leq \ro t + b \mid \lm \} = \Pr \li \{ \f{(\ro - \lm) t - a}{ \sq{\lm t} } \leq Y_t \leq \f{(\ro - \lm) t + b}{ \sq{\lm t} } \mid \lm
\ri \}.
\]
Noting that \bee &  & \bb{E} \li [ \exp \li (  s \f{ X_t - \lm t }{\sq{\lm t}} \ri )   \ri ] = \exp (- s \sq{\lm t}) \; \bb{E} \li [ \exp \li (
s \f{ X_t }{\sq{\lm t}} \ri )   \ri ]  = \exp \li \{ - s \sq{\lm t} + \lm t \li [ \exp \li ( \f{s}{\sq{\lm t}} \ri ) - 1 \ri ] \ri \}  \eee and
that \bee \lim_{t \to \iy} \li \{ - s \sq{\lm t} + \lm t \li [ \exp \li ( \f{s}{\sq{\lm t}} \ri ) - 1 \ri ] \ri \}  & =  & \lim_{t \to \iy}
\li \{ - s t + t^2 \li [ \exp \li ( \f{s}{t} \ri ) - 1 \ri ] \ri \}\\
&  = & \lim_{t \to \iy} \li \{  -s t + t^2 \li [ 1 + \f{s}{t} + \f{s^2}{2 t^2} + O \li ( \f{1}{t^3} \ri ) - 1 \ri ]  \ri \} = \f{s^2}{2}, \eee
we have that $Y_t = \f{ X_t - \lm t }{\sq{\lm t}}$ converges to a Gaussian random variable with zero mean and unit variance as $t \to \iy$.
Consequently,
\[
\Pr \{ \ro t - a \leq X_t \leq \ro t + b \mid \lm \} = \Pr \li \{ \f{(\ro - \lm) t - a}{ \sq{\lm t} } \leq Y_t \leq \f{(\ro - \lm) t + b}{
\sq{\lm t} } \mid \lm \ri \} \to 0
\]
as $t \to \iy$, which can be readily shown by considering the cases of $\lm < \ro, \; \lm = \ro$ and $\lm > \ro$ as in the proof of Theorem
\ref{genST} in Appendix \ref{genST_app}.

\epf

\beL

\la{suffPOs}

For arbitrary integer $n$ and real numbers $t_i, \; i = 0, \cd, n$ with $0 = t_0 < t_1 < \cd < t_{n-1} < t_n = t$, the conditional probability
mass function of $X_{t_i}, \; i = 0, 1, \cd, n - 1$ given the value of $X_t$ does not depend on $\lm$.

\eeL

\bpf

Note that for a Poisson process $(X_t)_{t \in [0, \iy)}$ with an arrival rate $\lm > 0$, we have \bee \f{ \Pr \{ X_{t_i} = x_i, \; i = 1, \cd, n
\} }{ \Pr \{ X_{t_n} = x_n \} } = \f{ \prod_{i=1}^n  \f{ [(t_i - t_{i-1})\lm]^{x_i - x_{i-1}} e^{-\lm (t_i - t_{i-1})} }{ (x_i - x_{i-1})!  } }
{ \f{ (t_n \lm)^{x_n} e^{-\lm t_n} }{x_n !}  } = \f{x_n !} { (t_n)^{x_n} } \prod_{i=1}^n  \f{ (t_i - t_{i-1})^{x_i - x_{i-1}} }{ (x_i -
x_{i-1})! }, \eee where $x_0 = 0$. This implies that the conditional PMF of $X_{t_i}, \; i = 1, \cd, n - 1$ given the value of $X_{t_n}$ does
not involve $\lm$.

\epf

We are now in a position to prove the theorem. It can be readily checked that $f_t (X_t; \lm)$ is unimodal with respect to $\lm > 0$. Applying
this fact and Lemma \ref{PoSTL} leads to the conclusion that the observational process will eventually terminate with probability $1$. As a
consequence of the proven termination property and Lemma \ref{suffPOs}, statements (I), (II) and (III) of Theorem \ref{PoST} follow from Theorem
\ref{MSPRT_Composite_Continuous}. This completes the proof of the theorem.

\sect{Proof of Theorem \ref{norST} }  \la{norST_app}

We need some preliminary results.

\beL \la{norST} For arbitrary $\al, \ba \in (0, 1)$ and $\mu, \mu^{\prime}, \mu^{\prime \prime} \in (-\iy, \iy)$ with $\mu^{\prime} <
\mu^{\prime \prime}$, {\small \bee \lim_{t \to \iy} \Pr \li \{ \ba \leq \Up_t (X_t; \mu^{\prime}, \mu^{\prime \prime}) \leq \f{1}{\al}  \mid \mu
\ri \} = 0. \eee}   \eeL

\bpf

Note that \[ \Pr \li \{ \ba \leq \Up_t (X_t; \mu^{\prime}, \mu^{\prime \prime}) \leq \f{1}{\al}  \mid \mu \ri \} = \Pr \li \{ \f{(\mu^{\prime} +
\mu^{\prime \prime}) t}{2} + \f{\si^2}{\mu^{\prime \prime} - \mu^{\prime}} \ln \ba \leq X_t \leq \f{(\mu^{\prime} + \mu^{\prime \prime}) t}{2} +
\f{\si^2}{\mu^{\prime \prime} - \mu^{\prime}} \ln \f{1}{\al} \mid \mu \ri \}.
\]
Thus, $\Pr \li \{ \ba \leq \Up_t (X_t; \mu^{\prime}, \mu^{\prime \prime}) \leq \f{1}{\al}  \mid \mu \ri \}$ can be written as $\Pr \{ \ro t - a
\leq X_t \leq \ro t + b \mid \mu \}$, where $\ro = \f{(\mu^{\prime} + \mu^{\prime \prime})}{2}$ and $a, \; b$ are some positive numbers.  Define
$Y_t = \f{ X_t - \mu t }{ \si \sq{ t}}$.  Then, $Y_t$ is a Gaussian random variable with zero mean and unit variance.  It follows that
\[
\Pr \{ \ro t - a \leq X_t \leq \ro t + b \mid \mu \} = \Pr \li \{ \f{(\ro - \mu)t - a}{ \si \sq{ t} } \leq Y_t \leq \f{(\ro - \mu)t + b}{ \si
\sq{ t} } \mid \mu \ri \} \to 0
\]
as $t \to \iy$,  which can be readily shown by considering the cases of  $\mu < \ro, \; \mu = \ro$ and $\mu > \ro$.

\epf

\beL

\la{sufnor}  For arbitrary integer $n$ and real numbers $t_i, \; i = 0, \cd, n$ with $0 = t_0 < t_1 < \cd < t_{n-1} < t_n = t$, the conditional
probability density function of $X_{t_i}, \; i = 0, 1, \cd, n - 1$ given the value of $X_t$ does not depend on $\mu$.  \eeL

\bpf

Define $Z_i = X_{t_i}  - X_{t_{i-1}}$ for $i = 1, \cd, n$, where $X_{t_0} = X_0 = 0$.  Then, $Z_i$ are independent Gaussian variables with PDFs
\bee f_i (z_i) = \f{1}{ \sq{2 \pi (t_i - t_{i-1}) } \si  } \exp \li ( - \f{[z_i - (t_i - t_{i-1}) \mu ]^2}{2 (t_i - t_{i-1}) \si^2}
 \ri ), \qu i = 1, \cd, n.   \eee
 Note that $\Pr \{ X_{t_i} \leq x_i, \; i = 1, \cd, n \} = \int \cd \int_{ (z_1, \cd, z_n) \in S}  \; \prod_{i = 1}^n f_i (z_i) \; dz_1 \cd
dz_n$,  where $S = \{ (z_1, \cd, z_n): \sum_{i = 1}^j z_i \leq x_j \; \tx{for} \; j = 1, \cd, n \}$.  Define $y_j = \sum_{i = 1}^j z_i$ for $j =
1, \cd, n$. Then, $z_1 = y_1$ and $z_j = y_j - y_{j-1}$ for $j = 2, \cd, n$. Note that the determinant of the Jacobian of the transformation is
equal to $1$ and thus \[ \Pr \{ X_{t_i} \leq x_i, \; i = 1, \cd, n \}  = \int_{-\iy}^{x_1} f_1(y_1)   \cd \int_{-\iy}^{x_j } f_j (y_j - y_{j-1})
\cd \int_{-\iy}^{x_n} f_n (y_n - y_{n-1} ) dy_n \cd dy_j \cd dy_1. \]  Sequentially taking partial derivatives of the multiple integral with
respect to $x_n, x_{n-1}, \cd, x_1$ gives
\[
\f{\pa^n }{\pa x_1 \pa x_2 \cd \pa x_n} \Pr \{ X_{t_i} \leq x_i, \; i = 1, \cd, n \} = \prod_{i=1}^n f_i (x_i - x_{i-1}), \qu x_0 \DEF 0.
\]
It can be checked that
\[
\f{ \f{\pa^n }{\pa x_1 \pa x_2 \cd \pa x_n} \Pr \{ X_{t_i} \leq x_i, \; i = 1, \cd, n \}  } { \f{\pa }{\pa x_n} \Pr \{ X_{t_n} \leq x_n \}  } =
\f{ \prod_{i=1}^n  \f{1}{ \sq{2 \pi (t_i - t_{i-1}) } \si  } \exp \li ( - \f{(x_i - x_{i-1})^2}{2 (t_i - t_{i-1}) \si^2}
 \ri ) } {  \f{1}{ \sq{2 \pi t_n } \si  } \exp \li ( - \f{x_n^2}{2 t_n  \si^2}
 \ri ) },
\]
which is independent of $\mu$.

\epf

We are now in a position to prove the theorem.  It can be readily checked that $f_t (X_t; \mu, \si)$ is unimodal with respect to $\mu \in (-
\iy, \iy)$.  This fact together with Lemma \ref{norST} lead to the conclusion that the observational process will eventually terminate with
probability $1$. As a consequence of the proven termination property and Lemma \ref{sufnor}, statements (I), (II) and (III) of Theorem
\ref{norST} follow from Theorem \ref{MSPRT_Composite_Continuous}. This completes the proof of the theorem.

\end{document}